\documentclass[12pt,a4paper]{article}

\unitlength\textwidth
\divide\unitlength by 300\relax

\usepackage{amsmath,amssymb}
\usepackage{bm}
\bmdefine{\NNN}{N}
\bmdefine{\ZZZ}{Z}
\bmdefine{\RRR}{R}
\bmdefine{\sss}{s}

\newcommand{\RRRRR}{{\mathcal R}}
\newcommand{\TTTTT}{{\mathcal T}}
\newcommand{\DDDDD}{{\mathcal D}}

\newcommand{\covers}{\mathrel{\cdot\!\!\!>}}
\newcommand{\covered}{\mathrel{<\!\!\!\cdot}}
\newcommand{\ini}{\mathop{\rm in}\nolimits}

\newcommand{\height}{\mathrm{ht}}

\newcommand{\define}{\mathrel{:=}}
\newcommand{\definebycond}{\stackrel{\mathrm{def}}{\iff}}

\newcommand{\gor}{Gorenstein}
\newcommand{\cm}{Cohen-Macaulay}

\newcommand{\qed}{\nolinebreak\rule{.3em}{.6em}}
\newcommand{\joinirred}{join-irreducible}
\newcommand{\rank}{\mathrm{rank}}
\newcommand{\meet}{\wedge}
\newcommand{\bigmeet}{\bigwedge}
\newcommand{\join}{\vee}

\newcommand{\lm}{{\rm lm}}

\newcommand{\relint}{\mathrm{relint}}
\newcommand{\msize}{\mathrm{size}}
\newcommand{\orth}{{\mathrm{O}}}
\newcommand{\sorth}{{\mathrm{SO}}}
\newcommand{\GL}{{\mathrm{GL}}}
\newcommand{\glin}{{\mathrm{GL}}}

\newcommand{\chara}{{\mathrm{char}}}
\newcommand{\hilb}{{\mathrm{Hilb}}}
\newcommand{\image}{{\mathrm{Im}}}
\newcommand{\transpose}{^\top}

\newtheorem{thm}{Theorem}[section]

\newtheorem{example}[thm]{Example}
\newtheorem{lemma}[thm]{Lemma}
\newtheorem{cor}[thm]{Corollary}
\newtheorem{definition}[thm]{Definition}
\newtheorem{prop}[thm]{Proposition}
\newtheorem{remark}[thm]{Remark}

\newcommand{\refeq}[1]{(\ref{#1})}
\numberwithin{equation}{section}

\newcounter{cond}
\numberwithin{cond}{section}
\renewcommand{\thecond}{(\arabic{section}.\arabic{cond})}

\newtheorem{case}{Case}[thm]

\def\mysubsection{}

\makeatletter
\newcommand{\mysloppy}{\tolerance 9999 \hfuzz .5\p@ \vfuzz .5\p@}
\makeatother
\title{%
Doset Hibi rings with an application to invariant theory}
\author{%
Mitsuhiro MIYAZAKI}
\date{%
Dept.\ Math.,
Kyoto University of Education,
\\
Fushimi-ku, Kyoto, 612-8522, Japan}

\begin{document}
\sloppy

\maketitle

\begin{abstract}
We define the concept of a doset Hibi ring and a generalized doset Hibi ring
which are subrings of a Hibi ring and are normal affine semigrouprings.
We apply the theory of (generalized) doset Hibi rings to analyze
the rings of absolute orthogonal invariants and absolute special orthogonal invariants
and show that these rings are normal and \cm\
and has rational singularities if the characteristic of the base field is zero
and is
F-rational otherwise.
We also state criteria of \gor\ property of these rings.
\\
Keywords:Hibi ring, sagbi, doset, ring of invariants,
Schubert subvariety
\\
MSC:13F50, 13A50, 13H10, 14M15
\end{abstract}

\section{Introduction}

Grassmannians and their Schubert subvarieties are fascinating objects
and attract many mathematicians.
Let $m$ and $n$ be integers with $1\leq m\leq n$ and
$G_{m}(V)$ be the Grassmannian consisting of all $m$-dimensional
vector subspaces of an $n$-dimensional vector space $V$ over a field $K$.
Then the homogeneous coordinate ring of $G_{m}(V)$ is $K[\Gamma(X)]$,
where $X=(X_{ij})$ is the $m\times n$ matrix of indeterminates
and $\Gamma(X)$ is the set of $m$-minors of $X$.

Fix a complete flag
$0=V_0\subsetneq V_1\subsetneq\cdots\subsetneq V_n=V$
of subspaces of $V$.

\begin{definition}\rm
For integers $a_1$, \ldots, $a_m$ with
$1\leq a_1<\cdots<a_m\leq n$, 
we define
$
\Omega(a_1,\ldots,a_m)\define
\{W\in G_m(V)\mid
\dim(W\cap V_{a_i})\geq i\mbox{ for $i=1$, $2$, \ldots, $m$}
\}$.
\end{definition}
It is known that 
$\Omega(a_1,\ldots,a_m)$ is a
subvariety of $G_m(V)$ and 
called the Schubert
subvariety of $G_m(V)$.

Set $b_i\define n-a_{m-i+1}+1$ for $i=1$, \ldots, $m$.
Then it is known that there exists a universal $m\times n$ matrix
$Z$ with entries in a $K$-algebra $S$ with the condition
$I_i(Z_{\leq b_i-1})=(0)$ for $i=1$, \ldots, $m$,
where $Z_{\leq j}$ denotes  the $m\times j$ matrix consisting of the first $j$ columns of $Z$
and $I_i(M)$ denotes the ideal generated by the $i$-minors of a matrix $M$.
To be pricise, 
$Z$ is an $m\times n$ matrix with entries in a $K$-algebra $S$ such that
\begin{enumerate}
\refstepcounter{cond}
\item[\thecond]
\label{cond:zero}
$I_i(Z_{\leq b_i-1})=(0)$ for $i=1$, \ldots, $m$ and
\refstepcounter{cond}
\item[\thecond]
\label{cond:univ}
if $I_i(M_{\leq b_i-1})=(0)$ for $i=1$, \ldots, $m$
for a matrix with entries in some $K$-algebra $S'$,
then there is a unique $K$-algebra homomorphism
$K[Z]\to S'$ which
maps $Z$ to $M$,
\end{enumerate}
where $K[Z]$ denotes the $K$-subalgebra of $S$ generated by the entries of $Z$.
By the universal property, $Z$ is unique up to isomorphism.
Further, it is also known that $K[\Gamma(Z)]$ is the homogeneous coordinate ring
of $\Omega(a_1,\ldots,a_m)$, where $\Gamma(Z)$ is the set of $m$-minors of $Z$.

On the other hand,
because of the universal property of $Z$, any subgroup of 
$\GL(m,K)$ acts on $K[Z]$.
In this paper, we study the rings of absolute 
$\orth(m)$ and $\sorth(m)$-invariants (see Definition \ref{def:inv}),
denoted as $K[Z]^{\orth(m,-)}$ and $K[Z]^{\sorth(m,-)}$ respectively,
of $K[Z]$.

We now state the contents of this paper.
We first establish notation,
recall known facts and state basic facts in Section \ref{sec:prel}
in order to simplify the proofs and clarify the structure of the theory
of the main part of this paper,
Sections \ref{sec:inv}, \ref{sec:doshibi} and \ref{sec:gor}.

In Section \ref{sec:inv}, we show that the ring of 
absolute 
$\orth(m)$ and $\sorth(m)$-invariants
are precisely the expected ones,
i.e., 
$K[Z]^{\orth(m,-)}=K[Z\transpose Z]$ and $K[Z]^{\sorth(m,-)}=K[Z\transpose Z,\Gamma(Z)]$.
See Subsection \ref{subsec:det ring} for notation.
We use the result of DeConcini and Procesi \cite{dp} about the ring of absolute
invariants.
Since Richman wrote in the introduction of  \cite{ric} that the proof of
\cite{dp} is incorrect, we state a down to earth proof of the result of DeConcini and Procesi
about the ring of absolute invariants.

In Section \ref{sec:doshibi}, 
we introduce the concept of doset Hibi rings and their generalization.
They are subrings of a Hibi ring and are normal affine semigrouprings,
therefore are \cm.
We apply the theory of
doset Hibi rings and generalized doset Hibi rings
to investigate the property of 
$K[Z]^{\orth(m,-)}$ (resp.\ $K[Z]^{\sorth(m,-)}$)
by showing it
is isomorphic to a subring of a polynomial ring over $K$,
whose initial algebra is a doset Hibi ring
(resp.\ a generalized doset Hibi ring)
in an appropriate monomial order.
In particular, $K[Z]^{\orth(m,-)}$ (resp.\  $K[Z]^{\sorth(m,-)}$)
is a normal \cm\ ring and 
has rational singularities if $\chara K=0$
and is
F-rational if $\chara K>0$.

Finally in Section \ref{sec:gor}, we state a criterion of \gor\ property of
a doset Hibi ring and a generalized doset Hibi ring.
As an application, we state a criterion of \gor\
preperty of $K[Z]^{\orth(m,-)}$ (resp.\ $K[Z]^{\sorth(m,-)}$),
see Theorems \ref{thm:gor orth} and \ref{thm:gor sorth}.
We note that the criterion of \gor\ property of
$K[Z]^{\orth(m,-)}$ is already obtained by Conca \cite{con}
and Goto \cite{got2}.
Therefore, Theorem \ref{thm:gor orth} 
recovers
their result.

Here we remark that 
for an $n\times n$ matrix of indeterminates
$T=(T_{ij})$,
Bruns, R\"omer and Wiebe \cite[Remark 3.10]{brw} 
pointed out that $K[T]/I_{m+1}(T)$ is isomorphic to $K[X\transpose X]$ and normal
\cm\ with rational singularities if $\chara K=0$ and is F-rational if 
$\chara K>0$
by using the same method as this paper.
Generalization of their result to $K[Z\transpose Z]$ and the result of 
$K[Z\transpose Z, \Gamma(Z)]$ are contained in Sections \ref{sec:prel} and \ref{sec:doshibi}.
We also give a criterion of the \gor\ property for these generalized cases in Section \ref{sec:gor}
and
study the relation to the invariant theory in Section \ref{sec:inv}.

\section{Preliminaries}
\label{sec:prel}

In this section, we establish notation, recall known facts and show some
basic facts in order to prepare the main part of this paper.

In this paper, all rings and algebras are commutative 
with identity element.
We denote by $\NNN$ the set of non-negative integers,
by $\ZZZ$ the set of integers,
by $\RRR$ the set of real numbers
and by $\RRR_{\geq0}$ the set of non-negative real numbers.
We denote by $|X|$ the cardinality of a set $X$.

\subsection{Normal affine semigrouprings}

\mysubsection{}
Let $K$ be a field and $X_1$, \ldots, $X_r$ indeterminates,
$S$ a finitely generated  additive submonoid of $\NNN^r$.
We set $K[S]\define K[X^\sss\mid\sss\in S]$ where $X^\sss=X_1^{s_1}\cdots X_r^{s_r}$
for $\sss=(s_1,\ldots, s_r)$.
Then
\begin{thm}[\cite{hoc}]
\label{thm:hoc normal}
\begin{enumerate}
\item
$K[S]$ is normal if and only if $S=\ZZZ S\cap\RRR_{\geq0}S$.
\item
If $K[S]$ is normal, then it is \cm.
\end{enumerate}
\end{thm}
If $S=\ZZZ S\cap\RRR_{\geq 0}S$, then by Farkas' theorem
and the fundamental theorem of abelian groups,
$S$ is isomorphic to a toroidal
monoid in the terminology of Stanley \cite[p.\ 81]{sta2}.
Therefore by \cite[p.\ 82]{sta2}, 
we have the follwoing
\begin{thm}\label{thm:sta can}
If $K[S]$ is normal, then the canonical module of $K[S]$ is
\[
\bigoplus_{\sss\in\relint \RRR_{\geq0}S\cap S} K X^\sss,
\]
where $\relint\RRR_{\geq0}S$ denotes for the interior of $\RRR_{\geq0}S$
in the affine subspace spanned by $\RRR_{\geq0}S$.
\end{thm}
Stanley also stated, in the same article, 
a criterion of \gor\ property for graded \cm\ domains.
\begin{thm}[{\cite[4.4 Theorem]{sta2}}]\label{thm:sta gor cri}
Let $A$ be an $\NNN$-graded \cm\ domain with $A_0=K$.
Suppose $\dim A=d$.
Then $A$ is \gor\ if and only if there is $\rho\in\ZZZ$ such that
$\hilb(A,1/\lambda)=(-1)^d\lambda^\rho\hilb(A,\lambda)$,
where $\hilb$ denotes the Hilbert series.
\end{thm}
As a corollary, we have the following fact.
\begin{cor}\label{cor:same hilb}
Let $A$ and $B$ be $\NNN$-graded \cm\ domains with $A_0=B_0=K$.
Suppose $\hilb(A,\lambda)=\hilb(B,\lambda)$.
Then $A$ is \gor\ if and only if $B$ is \gor.
In particular, if $A$ is a graded subring of a polynomial ring with monomial order
and the initial algebra $\ini A$ of $A$ is finitely generated over $K$ and \cm,
then $A$ is \gor\ if and only if $\ini A$ is \gor.
\end{cor}
Note, by the flat deformation argument, 
we see that if $A$ is an $\NNN$-graded subalgebra of a polynomial ring with monomial order
and $\ini A$ is \cm, 
then
$A$ is \cm.
%
%

\subsection{Determinantal rings}

\label{subsec:det ring}

\mysubsection{}
\label{mysubsec:mat}
Let $m$ and $n$ be integers
with $0<m\leq n$
and $K$ a field.
For an $m\times n$ matrix $M$ with entries in a
$K$ algebra $S$,
we denote by
$M\transpose $ the transposed matrix of $M$, 
by $I_t(M)$ the ideal of $S$ generated by  $t$-minors of $M$,
by $M_{\leq j}$ the  $m\times j$ matrix consisting of the first
$j$-columns of $M$,
by 
$\Gamma(M)$ the set of $m$-minors of $M$
and by
$K[M]$ the $K$-subalgebra of $S$ generated by the entries of $M$.

We define posets
$\Gamma(m\times n)$ and
$\Gamma'(m\times n)$ by
$\Gamma(m\times n)\define
\{[c_1,\ldots,c_m]\mid 1\leq c_1<\cdots<c_m\leq n\}$ and
$\Gamma'(m\times n)\define
\{[c_1,\ldots,c_r]\mid r\leq m, 1\leq c_1<\cdots<c_r\leq n\}$.
The order of 
$\Gamma'(m\times n)$ is defined by
\[
\begin{array}{rl}
&[c_1,\ldots,c_r]\leq[d_1,\ldots, d_s]\\
\stackrel{\rm def}{\Longleftrightarrow}&
r\geq s,\ c_i\leq d_i\mbox{ for $i=1$, $2$, \ldots, $s$,}
\end{array}
\]
and the order of 
$\Gamma(m\times n)$ is defined by that of 
$\Gamma'(m\times n)$.
Note that 
$\Gamma'(m\times n)$ is a distributive lattice and
$\Gamma(m\times n)$ is a sublattice of $\Gamma'(m\times n)$.
For $[c_1,\ldots, c_r]\in\Gamma'(m\times n)$,
we define its size to be $r$ and denote it by
$\msize[c_1,\ldots, c_r]$.

We also define the poset $\Delta(m\times n)$ by
$\Delta(m\times n)\define\{[\alpha|\beta]\mid
\alpha\in\Gamma'(m\times m),\beta
\in\Gamma'(m \times n), \msize\alpha=\msize\beta\}$,
and define the order of $\Delta(m\times n)$ 
by
\[
\begin{array}{rl}
&[\alpha|\beta]\leq[\alpha'|\beta']\\
\stackrel{\rm def}{\Longleftrightarrow}&
\alpha\leq\alpha'\mbox{ in $\Gamma'(m\times m)$ and }
\beta\leq\beta'\mbox{ in $\Gamma'(m\times n)$.}
\end{array}
\]
Note also that $\Delta(m\times n)$ is a 
distributive lattice.

For $\gamma\in\Gamma(m\times n)$, we set
$\Gamma(m\times n;\gamma)\define\{\delta\in\Gamma(m\times n)\mid
\delta\geq\gamma\}$.
$\Gamma'(m\times n;\gamma)$ for $\gamma\in\Gamma'(m\times n)$
and $\Delta(m\times n;\delta)$ for $\delta\in\Delta(m\times n)$
are defined similarly.

For an $m\times n$ matrix $M=(m_{ij})$ and
$\delta=[c_1,\ldots,c_r|d_1,\ldots, d_r]\in
\Delta(m\times n)$, we denote the minor
$\det(m_{c_id_j})$ by $\delta_M$ or $[c_1,\ldots,c_r|d_1,\ldots, d_r]_M$.
We also denote the maximal minor
$\det(m_{ic_j})$ of $M$
by $\delta_M$ or $[c_1,\ldots, c_m]_M$,
for
$\delta=[c_1,\ldots, c_m]\in\Gamma(m\times n)$.

We use standard terminology on Hodge algebras.
Standard references on the notion of  Hodge algebra are
\cite{dep2} and \cite[Chapter 7]{bh}.
However, we use the term
``algebra with straightening law'' (ASL for short)
to mean an ordinal Hodge algebra in the terminology of \cite{dep2}.

Let $X=(X_{ij})$ be an $m\times n$ matrix of indeterminates,
i.e.,
entries $X_{ij}$ are independent indeterminates.
Then the following facts are known.
\begin{thm}[\cite{dep1}]
\label{thm:dep1}
\begin{enumerate}
\item
$K[\Gamma(X)]$ is an ASL over $K$ on $\Gamma(m\times n)$ 
with structure map $\delta\mapsto\delta_X$.
\item
$K[X]$ is an ASL over $K$ on $\Delta(m\times n)$ with structure map
$[\alpha|\beta]\mapsto[\alpha|\beta]_X$.
\end{enumerate}
\end{thm}

From now on, 
we fix $\gamma=[b_1, b_2,\ldots, b_m]\in\Gamma(m\times n)$.
Set $\gamma^*=[1,\ldots, m| b_1,\ldots, b_m]\in\Delta(m\times n)$ and
$R=K[X]/(\Delta(m\times n)\setminus\Delta(m\times n;\gamma^*))K[X]$.
Then $R$ is an ASL on $\Delta(m\times n;\gamma^*)$.
Therefore, the pair $K$-algebra $R$ and the image of $X$ in $R$ satisfy
\ref{cond:zero} and \ref{cond:univ} of introduction.


\mysubsection{}
Set
\[
W\define
\mbox{\small$
\left(
\begin{array}{cccc}
W_{11}&W_{12}&\cdots&W_{1m}\\
W_{21}&W_{22}&\cdots&W_{2m}\\
\vdots&\vdots&\ddots&\vdots\\
W_{m1}&W_{m2}&\cdots&W_{mm}
\end{array}
\right)$}
\]
and
\[
U_\gamma\define
\mbox{\small$
\left(
\begin{array}{cccccccccccc}
0&\cdots&0&U_{1b_1}&\cdots&U_{1b_2-1}&U_{1b_2}&\cdots&\cdots&
U_{1b_m}&
\cdots&U_{1n}\\
0&\cdots&0&0&\cdots&0&U_{2b_2}&\cdots&\cdots&
U_{2b_m}&
\cdots&U_{2n}\\
\cdots&\cdots&\cdots&\cdots&\cdots&\cdots&\cdots&\cdots&\cdots&\cdots&\cdots\\
0&\cdots&0&0&\cdots&0&0&\cdots&0&U_{mb_m}&\cdots&U_{mn}
\end{array}
\right)$},
\]
where
$W_{ij}$  and  $U_{ij}$ are independent indeterminates.
Then
\begin{thm}[\cite{miy1}]
\label{thm:wugamma univ}
The natural map
$K[X]/(\Delta(m\times n)\setminus\Delta(m\times n;\gamma^*))K[X]\to K[WU_\gamma]$
induced by 
$K[X]\to K[WU_\gamma]$, ($X\mapsto WU_\gamma$)
is an isomorphism. 
In particular,
$WU_\gamma$ has the universal property 
\ref{cond:zero} and \ref{cond:univ} of introduction.
\end{thm}

Here we state a result on the initial algebra of $K[WU_\gamma]$ in $K[W,U_\gamma]$
along the same line as \cite[section 3]{brw} (see also \cite{bc}).
As a by-product, we obtain another proof of Theorem \ref{thm:wugamma univ}.
%
%
%
%

First we introduce  a diagonal monomial order
(i.e., a monomial order such that the leading monomial of any minor of $W$ or
$U_\gamma$ is the product of the entries of its main diagonal)
on $K[W, U_\gamma]$ 
with $W_{ij}>W_{\imath'j}$,
$W_{ij}>W_{i\jmath'}$,
$U_{ij}>U_{\imath'j}$ and
$U_{ij}>U_{i\jmath'}$
for any possible $i$, $\imath'$, $j$ and $\jmath'$ with $i<\imath'$ and $j<\jmath'$.
For example, 
degree lexicographic order on the polynomial ring
$K[W, U_\gamma]$
defined 
by
$W_{11}>W_{21}>\cdots>W_{m1}>W_{12}>\cdots>W_{mm}>
U_{1b_1}>U_{1b_1+1}>\cdots>U_{1n}>
U_{2b_2}>\cdots>U_{mn}$.
Then we obtain the following
\begin{lemma}\label{lem:lm of stand monom}
\begin{enumerate}
\item
\label{item:lm minor}
For
$[c_1,\ldots, c_r|d_1,\ldots, d_r]\in\Delta(m\times n;\gamma^*)$,
\[
\lm([c_1,\ldots,c_r|d_1,\ldots, d_r]_{WU_\gamma})
=\prod_{j=1}^r W_{c_j j}U_{j d_j}.
\]
\item
\label{item:lm muwu}
If
$\mu=\prod_{i=1}^u[c_{i1},\ldots, c_{ir(i)}|d_{i1},\ldots, d_{ir(i)}]$
is a standard monomial on $\Delta (m\times n;\gamma^*)$,
then
\[
\lm(\mu_{WU_\gamma})=
\prod_{i=1}^u\prod_{j=1}^{r(i)}W_{c_{ij}j}U_{jd_{ij}},
\]
where
$\mu_{WU_\gamma}\define
\prod_{i=1}^u[c_{i1},\ldots, c_{ir(i)}|d_{i1},\ldots, d_{ir(i)}]_{WU_\gamma}$.
\item
If $\mu$ and $\mu'$ are standard monomials on $\Delta(m\times n;\gamma^*)$ with $\mu\neq\mu'$,
then
$\lm(\mu_{WU_\gamma})\neq\lm(\mu'_{WU_\gamma})$.
In particular,
$\{\mu_{WU_\gamma}\mid\mu$ is a standard monomial on $\Delta(m\times n;\gamma^*)\}$
is linearly independent over $K$.
\item
\label{item:xwu isom}
The natural map
$K[X]/(\Delta(m\times n)\setminus\Delta(m\times n;\gamma^*))K[X]\to K[WU_\gamma]$
induced by 
$K[X]\to K[WU_\gamma]$, ($X\mapsto WU_\gamma$)
is an isomorphism. 
In particular, $K[WU_\gamma]$ is an ASL on $\Delta(m\times n;\gamma^*)$.
\item
\label{item:unique stand monom}
For any $a\in K[WU_\gamma]$,
there is a unique standard monomial $\mu$ on $\Delta(m\times n;\gamma^*)$
such that
$\lm(a)=\lm(\mu_{WU_\gamma})$.
\end{enumerate}
\end{lemma}
\begin{proof}
\leavevmode
\ref{item:lm minor} is proved 
by using  the
following equation on determinants and the assumption of the monomial order.
\begin{eqnarray*}
&&
[c_1,\ldots,c_r|d_1, \ldots, d_r]_{WU_\gamma}\\
&=&
\sum_{1\leq j_1<\cdots<j_r\leq m}
[c_1,\ldots,c_r|j_1, \ldots, j_r]_{W}
[j_1,\ldots,j_r|d_1, \ldots, d_r]_{U_\gamma}.
\end{eqnarray*}
\ref{item:lm muwu} follows from \ref{item:lm minor}.
The rest is proved 
straight forwardly 
or 
by the same way as \cite[Section 3]{brw}.
\end{proof}

We set $Z_\gamma\define WU_\gamma$ in the following.
As a corollary, we see the following

\begin{prop}
$Z_\gamma$ has the universal property \ref{cond:zero} and \ref{cond:univ} and
$\{[\alpha|\beta]_{Z_\gamma}\mid
[\alpha|\beta]\in\Delta(m\times n;\gamma^*)\}$ is
a sagbi basis of $K[Z_\gamma]$.
\end{prop}

\mysubsection{}
We set $D_{m,n}\define\{(\alpha,\beta)\mid\alpha$, $\beta\in\Gamma'(m\times n)$,
$\alpha\leq\beta$ and 
$\msize\alpha=\msize\beta\}$ and
$D'_{m,n}\define D_{m-1,n}\cup\Gamma(m\times n)$.
We define the
order on $D_{m,n}$ by lexcographic way, i.e.,
\[
(\alpha,\beta)<(\alpha',\beta')
\stackrel{\rm def}{\Longleftrightarrow}
\alpha<\alpha' \mbox{ or } (\alpha=\alpha'\mbox{ and }\beta<\beta').
\]
We also define the order on $D'_{m,n}$ by extending the orders of $D_{m-1,n}$
and $\Gamma(m\times n)$ and set
\[
\delta<(\alpha,\beta)\definebycond \delta<\alpha\mbox{ in $\Gamma'(m\times n)$}
\]
for $\delta\in\Gamma(m\times n)$ and $(\alpha,\beta)\in D_{m-1,n}$.
Then the following fact is known.
\begin{thm}[{\cite[Section 5]{dp}}]
\label{thm:dp sec5}
\begin{enumerate}
\item
$K[X\transpose X]$
is a Hodge algebra over $K$ on $D_{m,n}$ with structure map
$(\alpha,\beta)\mapsto[\alpha|\beta]_{X\transpose X}$,
where the ideal of monomials on $D_{m,n}$ defining this Hodge algebra structure is generated
by
\[
\{(\alpha,\beta)(\alpha',\beta')\mid \beta\not\leq \alpha'\mbox{ and }\beta'\not\leq\alpha\}.
\]
\item
$K[X\transpose X,\Gamma(X)]$
is a Hodge algebra over $K$ on $D'_{m,n}$ with structure map
$\delta\mapsto \delta_X$ for $\delta\in\Gamma(m\times n)$ and
$(\alpha,\beta)\mapsto[\alpha|\beta]_{X\transpose X}$ for $(\alpha,\beta)\in D_{m-1,n}$,
where the ideal of monomials on $D'_{m,n}$ defining this Hodge algebra structure is generated
by
\begin{eqnarray*}
&&\{(\alpha,\beta)(\alpha',\beta')\mid 
(\alpha,\beta),(\alpha',\beta')\in D_{m-1,n},
\beta\not\leq \alpha'\mbox{ and }\beta'\not\leq\alpha\}\\
&\cup&
\{\delta(\alpha,\beta)\mid\delta\in\Gamma(m\times n), (\alpha,\beta)\in D_{m-1,n}\mbox{ and }
\delta\not\leq\alpha\mbox{ in }\Gamma'(m\times n)\}\\
&\cup&
\{\delta\delta'\mid\delta,\delta'\in\Gamma(m\times n), \delta\not\leq\delta',
\delta'\not\leq\delta\}.
\end{eqnarray*}
\end{enumerate}
\end{thm}

In the following, when we mention about the concept of Hodge algebra 
on $D_{m,n}$ and $D'_{m,n}$, such as standard monomials,
the ideals of monomials defining the Hodge algebra structures are the ones 
in Theorem \ref{thm:dp sec5}.
In \cite[Section 18]{dep2}, 
DeConcini, Eisenbud and Procesi
showed that $K[X\transpose  X]$ is a doset algebra.
However, 
every doset algebra has a structure of a Hodge algebra by ordering the
elements of the doset lexicographically.
%
%
In fact, 
DeConcini-Procesi \cite[Lemmas 5.2, 5.3 and 5.4]{dp}
essentially proved the following fact also.

\begin{thm}
\label{thm:dp sec5e}
Let $T=(T_{ij})$ be an $n\times n$ symmetric matrix of indeterminates,
i.e.,
$\{T_{ij}\}_{1\leq i\leq j\leq n}$ are independent indeterminates and
$T_{ji}=T_{ij}$ for $j>i$.
Then
$K[T]/I_{m+1}(T)$ is a Hodge algebra over $K$ on $D_{m,n}$ with structure map
$(\alpha,\beta)\mapsto[\alpha|\beta]_{\overline T}$,
where $\overline T$ is the image of $T$ in $K[T]/I_{m+1}(T)$
and we set $I_{m+1}(T)=(0)$ if $m=n$.
\end{thm}
As a corollary of Theorems \ref{thm:dp sec5} and \ref{thm:dp sec5e},
we see the following fact.

\begin{cor}
Let $T$ be an $n\times n$ symmetric matrix of indeterminates.
Then the $K$-algebra homomorphism
$K[T]/I_{m+1}(T)\to K[X]$ induced by $T\mapsto X\transpose X$ is an isomorphism.
\end{cor}

For an $m\times n$ matrix $M$ and for a standard monomial 
$\nu=\prod_{i=1}^v[c_{i1},\ldots,c_{ir(i)}],[d_{i1},\ldots,d_{ir(i)}]$
on $D_{m,n}$, we set
$\nu_M\define
\prod_{i=1}^v[c_{i1},\ldots,c_{ir(i)}|d_{i1},\ldots,d_{ir(i)}]_{M\transpose M}$
and for a standard monomial
$\nu'=\prod_{i=1}^{v'}[e_{i1},\ldots, e_{im}]
\prod_{i=1}^{v''}[c_{i1},\ldots,c_{ir(i)}],[d_{i1},\ldots,d_{ir(i)}]$
on $D'_{m,n}$, we set
$\nu'_M\define
\prod_{i=1}^{v'}[e_{i1},\ldots,e_{im}]_M
\prod_{i=1}^{v''}[c_{i1},\ldots,c_{ir(i)}|d_{i1},\ldots,d_{ir(i)}]_{M\transpose M}$.

By the same argument as the proof of Lemma~\ref{lem:lm of stand monom},
we see the following facts.
\begin{lemma}\label{lem:lm of symm}
\begin{enumerate}
\item
For $[c_1,\ldots, c_r|d_1,\ldots, d_r]\in\Delta(m\times m)$, we have
\[\lm([c_1,\ldots,c_r|d_1,\ldots,d_r]_{W\transpose W})
=\prod_{j=1}^r W_{jc_j}W_{jd_j}.
\]
\item
\label{item:lm of symm}
If
$\nu=\prod_{i=1}^v([c_{i1},\ldots,c_{ir(i)}],[d_{i1},\ldots,d_{ir(i)}])$
is a standard monomial on $D_{m,m}$, then
\[
\lm(\nu_{W})=
\prod_{i=1}^v\prod_{j=1}^{r(i)}W_{jc_{ij}}W_{jd_{ij}}.
\]
\item
If $\nu$ and $\nu'$ are standard monomials on $D_{m,m}$ with $\nu\neq\nu'$,
then $\lm(\nu_{W})\neq\lm(\nu'_{W})$.
In particular,
$\{\nu_{W}\mid\nu$ is a standard monomial on $D_{m,m}\}$
is linearly independent over $K$.
\item
\label{item:unique lm of symm}
For any $a\in K[W\transpose W]$, there is a unique standard monomial $\nu$ on $D_{m,m}$
such that
$\lm(a)=\lm(\nu_{W})$.
In particular, 
$\{[\alpha|\beta]_{W\transpose W}\mid(\alpha,\beta)\in D_{m,m}\}$
is a sagbi basis of $K[W\transpose W]$.
\end{enumerate}
\end{lemma}
\begin{lemma}\label{lem:lm of symm+}
\begin{enumerate}
\item
\label{item:lm of symm+}
If
$\nu=[1,\ldots,m]^{v'}
\prod_{i=1}^{v''}([c_{i1},\ldots,c_{ir(i)}],[d_{i1},\ldots,d_{ir(i)}])$
is a standard monomial on $D'_{m,m}$, then
\[
\lm(\nu_{W})=
(\prod_{j=1}^m W_{jj})^{v'}
\prod_{i=1}^{v''}\prod_{j=1}^{r(i)}W_{jc_{ij}}W_{jd_{ij}}.
\]
\item
If $\nu$ and $\nu'$ are standard monomials on $D'_{m,m}$ with $\nu\neq\nu'$,
then $\lm(\nu_{W})\neq\lm(\nu'_{W})$.
In particular,
$\{\nu_{W}\mid\nu$ is a standard monomial on $D'_{m,m}\}$
is linearly independent over $K$.
\item
\label{item:unique lm of symm+}
For any $a\in K[W\transpose W,\det W]$, there is a unique standard monomial $\nu$ on $D'_{m,m}$
such that
$\lm(a)=\lm(\nu_{W})$.
In particular, 
$\{\det W\}\cup\{[\alpha|\beta]_{W\transpose W}\mid(\alpha,\beta)\in D_{m-1,m}\}$
is a sagbi basis of $K[W\transpose W,\det W]$.
\end{enumerate}
\end{lemma}

Finally in this section we make the following remark.
Let $T$ be an $n\times n$ symmetric matrix of indeterminates.
Set
\[
R=K[T]/\{(\alpha,\beta)\in D_{n,n}\mid\alpha\not\geq\gamma\}K[T],
\]
and
$$
D_\gamma\define\{(\alpha,\beta)\in D_{n,n}\mid\alpha\geq\gamma\}.
$$
where $\gamma=[b_1, \ldots, b_m]\in\Gamma(m\times n)$ is the element
fixed after Theorem \ref{thm:dep1}.
Then $R$ is a Hodge algebra over $K$ on 
$D_\gamma$
since $D_{n,n}\setminus D_\gamma$ is a poset ideal of $D_{n,n}$.
Let $T_\gamma$ be the image of $T$ in $R$.
Then 
$T_\gamma$ has the following universal property,
where we set $b_{m+1}\define n+1$.
\begin{enumerate}
\item
$I_i((T_\gamma)_{\leq b_i-1})
=(0)$ for $i=1$, \ldots, $m+1$ and
\item
if $M$ is an $n\times n$ symmetric matrix with entries in some $K$-algebra $S$
such that
$I_i(M_{\leq b_i-1})
=(0)$ 
for $i=1$, \ldots, $m+1$,
then there is a unique $K$-algebra homomorphism
$K[T_\gamma]\longrightarrow S$ which
maps $T_\gamma$ to $M$.
\end{enumerate}

By the same argument as the proof of Lemma \ref{lem:lm of stand monom}, we see
the following fact,
cf.\ \cite[Remark 3.10]{brw}.
\begin{lemma}\label{lem:lm of zgamma t}
\begin{enumerate}
\item
For $([c_1,\ldots, c_r], [d_1,\ldots, d_r])\in D_\gamma$,
\[
\lm([c_1,\ldots,c_r|d_1,\ldots,d_r]_{Z_\gamma\transpose Z_\gamma})
=
\prod_{j=1}^r W_{jj}^2U_{jc_j}U_{jd_j}.
\]
\item
\label{item:lm of zgamma t}
If
$\nu=\prod_{i=1}^l([c_{i1},\ldots,c_{ir(i)}],[d_{i1},\ldots,d_{ir(i)}])$
is a standard monomial on $D_\gamma$, then
\[
\lm(\nu_{Z_\gamma})=\prod_{i=1}^l\prod_{j=1}^{r(i)}
W_{jj}^2U_{jc_{ij}}U_{jd_{ij}}.
\]
\item
If $\nu$ and $\nu'$ are standard monomials on $D_\gamma$ with $\nu\neq\nu'$, then
$\lm(\nu_{Z_\gamma})\neq\lm(\nu'_{Z_\gamma})$.
In particular,
$\{\nu_{Z_\gamma}\mid\nu$ is a standard monomial on $D_\gamma\}$ is 
linearly independent over $K$.
\item
\label{item:lm of zgamma t isom}
The natural map
$K[T]/(D_{n,n}\setminus D_\gamma)K[T]\to
K[Z_\gamma\transpose  Z_\gamma]$
induced by
$
K[T]\to K[Z_\gamma\transpose  Z_\gamma]$,
($T\mapsto Z_\gamma\transpose  Z_\gamma$)
is an isomorphism.
\item
$Z_\gamma\transpose  Z_\gamma$ is the universal $n\times n$ symmetric matrix with
$I_i(M_{\leq b_i-1})
=(0)$ for $i=1$, \ldots, $m+1$ and
$\{[\alpha|\beta]_{Z_\gamma\transpose  Z_\gamma}\mid (\alpha,\beta)\in D_\gamma\}$
is a sagbi basis of 
$K[Z_\gamma\transpose  Z_\gamma]
$.
\end{enumerate}
\end{lemma}

\section{Rings of invariants and sagbi bases}
\label{sec:inv}

Recall that we have fixed 
integers $m$, $n$ with $0<m\leq n$ and
$\gamma=[b_1,b_2,\ldots, b_m]\in\Gamma(m\times n)$.

For any $g\in\GL(m,K)$,
$I_i((gZ_\gamma)_{\leq b_i-1})=(0)$ for $i=1$, \ldots, $m$.
Therefore,
by the universal property of $Z_\gamma$,
 there is a $K$-automorphism of $K[Z_\gamma]$
sending $Z_\gamma$ to $gZ_\gamma$,
i.e., any subgroup of $\GL(m,K)$ acts on 
$K[Z_\gamma]$.

We define the orthogonal group and the special orthogonal group naively.

\begin{definition}
\rm
Let $R$ be a commutative ring.
We set $\orth(m,R)\define\{g\in \glin(m,R)\mid
g\transpose g=E_m\}$ and
$\sorth(m,R)\define\{g\in\orth(m,R)\mid\det g=1\}$.
\end{definition}

Next we make the following

\begin{definition}
\label{def:inv}
\rm
Let $f\in K[Z_\gamma]$.
$f$ is called a naive $\orth(m)$-invariant if $f^g=f$ for any
$g\in \orth(m,K)$.
$f$ is called an absolute $\orth(m)$-invariant if for any
$K$-algebra $B$, if we denote the image of $f$ in $B[Z_\gamma]$ as $f$,
$f^g=f$ for any $g\in\orth(m,B)$.
The ring of naive (resp.\ absolute) $\orth(m)$-invariants is denoted
by $K[Z_\gamma]^{\orth(m,K)}$ (resp.\ $K[Z_\gamma]^{\orth(m,-)}$).
Naive (resp.\ absolute) $\sorth(m)$-invarinats and $K[Z_\gamma]^{\sorth(m,K)}$
(resp.\ $K[Z_\gamma]^{\sorth(m,-)}$) are defined similarly.
\end{definition}
Note that we allow nilpotent elements in $B$ in the above definition.
In particular, $\sorth(m,B)$ is not identical to $\orth(m,B)$ 
even in the case where $\chara K=2$.
The concept of absolute $\orth(m)$ and $\sorth(m)$-invariants are 
defined in \cite[Section 2]{dp}.
We recall known results about the ring of naive and absolute 
$\orth(m)$ and $\sorth(m)$-invariants.

\begin{thm}
\label{thm:rdp}
Let $X=(X_{ij})$ be an $m\times n$ matrix of indeterminates.
\begin{enumerate}
\item
\label{item:non char 2}
If $K$ is an infinite field and $\chara K\neq 2$, then
$K[X]^{\orth(m,K)}=K[X\transpose  X]$
and
$K[X]^{\sorth(m,K)}=K[X\transpose  X,\Gamma(X)]$.
\item
\label{item:char 2}
If $K$ is an infinite field and $\chara K=2$, then
$K[X]^{\orth(m,K)}=K[X\transpose  X, 
\{X_{1j}+\cdots+X_{mj}\mid 1\leq j\leq n\}]$.
\item
\label{item:abs}
$K[X]^{\orth(m,-)}=K[X\transpose  X]$ and
$K[X]^{\sorth(m,-)}=K[X\transpose  X,\Gamma(X)]$.
\end{enumerate}
\end{thm}
\ref{item:non char 2} and \ref{item:char 2} are due to Richman \cite{ric}
and \ref{item:abs} is due to 
DeConcini and Procesi \cite{dp}.
\ref{item:abs} is valid even in the case where $\chara K=2$.
As is noted above, we require $f$ to be an absolute $\orth(m)$-invariant
that for any $K$-algebra $B$,
 including the one such that $\sorth(m,B)\subsetneq \orth(m,B)$,
$f^g=f$ for any $g\in \orth(m,B)$.
Thus, 
$K[X]^{\orth(m,-)}$ and $K[X]^{\sorth(m,-)}$ are different even in the
case where $\chara K=2$ 
(see the following another proof of \ref{item:abs} of Theorem \ref{thm:rdp}).
Since Richman \cite[page 45]{ric} wrote that the proof of \cite{dp} is incorrect,
 we state a down to earth proof of 
\ref{item:abs} of Theorem \ref{thm:rdp}.

First we show that $K[X]^{\sorth(m,-)}=K[X\transpose  X,\Gamma(X)]$.
We follow \cite[Section 4]{ric}.
Only the proof of \cite[Proposition 13]{ric} is to be modified.
We use the notation of \cite[Section 4]{ric}.
Suppose $f\in K[C, m\times 2]^{\sorth(2,-)}$.
Let $T_1$ and $T_2$ be indeterminates.
Set $B_1=K[T_1,T_2]/(T_1^2+1-T_2^2)$.
Since $T_2$ is a non-zerodivisor of $B_1$,
$B_1$ is a subring of $B_1[1/T_2]$.
We set $B=B_1[1/T_2]$
and 
$$
g=\frac{1}{T_2}\begin{pmatrix}1&-T_1\\T_1&1\end{pmatrix}.
$$
Then $g\in\sorth(m,B)$.

Since $f$ is an absolute $\sorth(2)$-invariant in the terminology of 
the present paper, we see that
$f^g=f$.
Write 
$$
f=c_1|M_1|\langle W_1\rangle+\cdots
+c_k|M_k|\langle W_k\rangle
$$
as \cite[(4.5)]{ric}.
Since
$$
(|M|\langle W\rangle)^g
=\left(\frac{T_1}{T_2}\right)^{e_1(M)}
\left(\frac{1}{T_2}\right)^{e_2(M)}
|\hat M|\langle W\rangle
+\sum_{{i\geq 0 \atop 0\geq j\geq -e_1(M)-e_2(M)}}h_{ij}T_1^i T_2^j
,
$$
where $h_{ij}\in K[C,m\times 2]$ and the leading monomials of
${h_{ij}}'s$ are strictly less than the leading monomial of
$|\hat M|\langle W\rangle$.
Thus,
if we define $E$, $J$ and $J'$ as in \cite[page 62]{ric},
$$
T_2^E f\equiv
\sum_{j\in J}c_j|\hat M_j|\langle W_j\rangle
+T_1\sum_{j\in J'}c_j|\hat M_j|\langle W_j\rangle
+y+T_1z
\mod (T_2,T_1^2+1),
$$
where $y$ and $z$ are elements of $K[C, m\times 2]$ whose leading
monomials are strictly less than the leading monomial of some element
$|\hat M_j|\langle W_j\rangle$ with $j\in J\cup J'$.

Since $1$, $T_1$ is a free basis of 
$K[C, m\times 2][T_1,T_2]/(T_2,T_1^2+1)$ as a 
free module over $K[C, m\times 2]$, we see that $E=0$.

Thus, we see the counterpart of \cite[Proposition 13]{ric},
and therefore, we see that
$K[X]^{\sorth(m,-)}=K[X\transpose X,\Gamma(X)]$
by the same argument as the proof of \cite[Proposition 14]{ric}.

Next we show that $K[X]^{\orth(m,-)}=K[X\transpose X]$.
The inclusion $K[X]^{\orth(m,-)}\supset K[X\transpose X]$ is obvious.
Suppose $f\in K[X]^{\orth(m,-)}$.
Then $f\in K[X]^{\sorth(m,-)}$.
Therefore, by the above proved fact and Theorem \ref{thm:dp sec5},
we see that there are standard monomials $\nu_1$, \ldots, $\nu_k$ of
$D'_{m,n}$ such that
$$
f=c_1(\nu_1)_X+\cdots+c_k(\nu_k)_X.
$$
Set $J=\{j\mid \nu_j$ is a product of odd number of elements of $\Gamma(m\times n)$
and elements of $D_{m-1,n}\}$, where we alow empty product of elements of 
$D_{m-1,n}$.
Let $T$ be an indeterminate and $B=K[T]/(T^2-1)$.
Then
$$
h=\begin{pmatrix}T\\&1&\\&&\ddots\\&&&1\end{pmatrix}\in\orth(m,B)
$$
and
$$
f^h=T\sum_{j\in J}c_j(\nu_j)_X+\sum_{j\not\in J}c_j(\nu_j)_X.
$$
Since $1$, $T$ is a free basis of $B[X]$ as a free module over $K[X]$ and
$f^h=f$, we see that $J=\emptyset$.
Therefore, $f\in K[X\transpose X]$.

Since the ring of naive invariants is  identical to the absolute invariant
except the case where $\chara K=2$ or $K$ is a finite field,
we focus on the absolute invariants in the rest of this paper.

As a corollary of \ref{item:abs} of Theorem \ref{thm:rdp}, we see the
following fact.

%
%
%
\begin{cor}\label{cor:inv of kw}
\begin{enumerate}
\item
\label{item:inv of kw orth}
$K[W]^{\orth(m,-)}=K[W\transpose W]$.
\item
$K[W]^{\sorth(m,-)}=K[W\transpose W,\det W]$.
\end{enumerate}
\end{cor}

Here we note the following basic fact.
\begin{lemma}
\label{lem:sagbi plus}
Let $K[X_1$, \ldots, $X_s$, $Y_1$, \ldots, $Y_t]$ be a polynominal
ring over $K$ with monomial order and $A$ a $K$-subalgebra of
$K[X_1$, \ldots, $X_s]$.
If $S$ is a sagbi basis of $A$ in $K[X_1$, \ldots, $X_s]$, then
$S\cup\{Y_1$, \ldots, $Y_t\}$ is a sagbi basis of 
$A[Y_1$, \ldots, $Y_t]$ in $K[X_1$, \ldots, $X_s$, $Y_1$, \ldots, $Y_t]$.
\end{lemma}
Now we state the following
\begin{thm}\label{thm:sagbi of inv}
\begin{enumerate}
\item
$\{[\alpha|\beta]_{Z_\gamma\transpose  Z_\gamma}\mid
(\alpha,\beta)\in D_\gamma\}$
is a sagbi basis of 
$K[W\transpose W,U_\gamma]\cap K[Z_\gamma]$.
In particular,
$K[Z_\gamma]^{\orth(m,-)}=
K[W\transpose W,U_\gamma]\cap K[Z_\gamma]=
K[Z_\gamma\transpose  Z_\gamma]$.

\item
$\{[\alpha|\beta]_{Z_\gamma\transpose  Z_\gamma}\mid
(\alpha,\beta)\in D_\gamma\}
%
\cup\{\delta_{Z_\gamma}\mid
\delta\in\Gamma(m\times n;\gamma)\}
$
is a sagbi basis of 
$K[W\transpose W, \det W,U_\gamma]\cap K[Z_\gamma]$.
In particular,
$K[Z_\gamma]^{\sorth(m,-)}=
K[W\transpose W,\det W, U_\gamma]\cap K[Z_\gamma]=
K[Z_\gamma\transpose  Z_\gamma, \Gamma(Z_\gamma)]$.
\end{enumerate}
\end{thm}
\begin{proof}
(1):
It is clear that
$K[W\transpose W, U_\gamma]\cap K[Z_\gamma]\supset K[Z_\gamma\transpose Z_\gamma]$.
Thus,  
$[\alpha|\beta]_{Z_\gamma\transpose Z_\gamma}\in
K[W\transpose W, U_\gamma]\cap K[Z_\gamma]$
for any $(\alpha,\beta)\in D_\gamma$.
Let $a$ be an arbitrary element of 
$K[W\transpose W, U_\gamma]\cap K[Z_\gamma]$.
Then by Lemma \ref{lem:lm of stand monom} \ref{item:unique stand monom}, we see that there is a
unique standard  monomial $\mu$ on $\Delta(m\times n;\gamma^*)$
 such that
$\lm(a)=\lm(\mu_{Z_\gamma})$.

Set
$\mu=\prod_{i=1}^u[\alpha_i|\beta_i]$,
$\alpha_1\leq\cdots\leq\alpha_u$,
$\beta_1\leq\cdots\leq\beta_u$,
and
$\alpha_i=[c_{i1},\ldots,c_{ir(i)}]$,
$\beta_i=[d_{1i},\ldots,d_{ir(i)}]$
for
$i=1$, \ldots, $u$.
Then
by Lemma \ref{lem:lm of stand monom} \ref{item:lm muwu}
we see that
\begin{equation}\label{eqn:asl zgamma}
\lm(a)=\lm(\mu_{Z_\gamma})
=\prod_{i=1}^u\prod_{j=1}^{r(i)}W_{c_{ij}j}U_{jd_{ij}}.
\end{equation}
On the other hand,
we see, by Lemma \ref{lem:lm of symm} \ref{item:lm of symm}, \ref{item:unique lm of symm}
and Lemma \ref{lem:sagbi plus}, that
\begin{equation}\label{eqn:hodge wtw}
\lm(a)=
\left(
\prod_{i=1}^v\prod_{j=1}^{s(i)}W_{jc'_{ij}}W_{jd'_{ij}}
\right)\times
\mbox{(a monomial of ${U_{\cdot\cdot}}'s$)}
\end{equation}
for some $c'_{ij}$ and $d'_{ij}$,
since $a\in K[W\transpose W, U_\gamma]$.

By comparing \refeq{eqn:asl zgamma} and \refeq{eqn:hodge wtw},
we see that
$c_{ij}=j$, $c'_{ij}=j$, and $d'_{ij}=j$ for any $i$, $j$ in the equations.
Further, we see that
$u$ is even and $r(2i-1)=r(2i)$ for $i=1$, \ldots, $u/2$.
Therefore if we set
$\alpha'_i=[d_{2i-1,1},\ldots, d_{2i-1,r(2i-1)}]$
and
$\beta'_i=[d_{2i,1},\ldots, d_{2i,r(2i)}]$,
then
$(\alpha'_i$, $\beta'_i)\in D_\gamma$ for $i=1$, \ldots, $u/2$,
$\prod_{i=1}^{u/2}(\alpha'_i,\beta'_i)$ is a standard monomial on
$D_{m,n}$
and
\[
\lm(a)=\prod_{i=1}^{u/2}\lm([\alpha'_i|\beta'_i]_{Z_\gamma\transpose Z_\gamma})
\]
by Lemma \ref{lem:lm of zgamma t} \ref{item:lm of zgamma t}.
Therefore, the former part of (1) follows.
Since the action of $\orth(m,B)$ on $B[Z_\gamma]$ is induced by the action of
$\orth(m,B)$ on $B[W]$ for any $K$-algebra $B$, we see the latter part of (1)
by Corollary \ref{cor:inv of kw} \ref{item:inv of kw orth}.

(2):
Again, it is clear that $K[W\transpose W, \det W, U_\gamma]\cap K[Z_\gamma]
\supset K[Z_\gamma\transpose Z_\gamma,\Gamma(Z_\gamma)]$.
Thus, $[\alpha|\beta]_{Z_\gamma}$, $\delta_{Z_\gamma}\in K[W\transpose W,\det W, U_\gamma]
\cap K[Z_\gamma]$ for any $(\alpha, \beta)\in D_\gamma$ and
$\delta\in \Gamma(m\times n;\gamma)$.
Let $a$ be an arbitrary element of 
$K[W\transpose W, \det W, U_\gamma]\cap K[Z_\gamma]$.
By 
Lemma \ref{lem:lm of stand monom} \ref{item:unique stand monom},
we see that there is a standard monomial
$\mu$
on $\Delta (m\times n;\gamma^*)$ such that
$\lm(a)=\lm(\mu_{Z_\gamma})$.
Set, as before,
$\mu=\prod_{i=1}^u[\alpha_i|\beta_i]$,
$\alpha_1\leq\cdots\leq\alpha_u$,
$\beta_1\leq\cdots\leq\beta_u$,
and
$\alpha_i=[c_{i1},\ldots,c_{ir(i)}]$,
$\beta_i=[d_{1i},\ldots,d_{ir(i)}]$
for
$i=1$, \ldots, $u$.
Then by 
using Lemma~\ref{lem:lm of symm+} \ref{item:lm of symm+}, \ref{item:unique lm of symm+},
Lemma~\ref{lem:sagbi plus}
and
the same argument 
as in the proof of (1), we see that
$c_{ij}=j$ for any $i$, $j$ and
 there is $w\in\NNN$
such that
$w\leq u$,
$r(i)=m\iff i\leq w$,
$u-w$ is even
and
$r(2i-1+w)=r(2i+w)$ for $i=1$, \ldots, $(u-w)/2$.
Therefore if we set
$\alpha''_i=[d_{2i-1+w,1},\ldots,d_{2i-1+w,r(2i-1+w)}]$
and
$\beta''_i=[d_{2i+w,1},\ldots,d_{2i+w,r(2i+w)}]$
for $i=1$, \ldots, $(u-w)/2$ and
$\delta_i=[d_{i1},\ldots, d_{im}]$ for $i=1$, \ldots, $w$,
then
\[
\lm(a)=
\prod_{i=1}^{u}\prod_{j=1}^{r(i)}W_{jj}U_{jd_{ij}}
=
\prod_{i=1}^w \lm((\delta_i)_{Z_\gamma})
\prod_{i=1}^{(u-w)/2} \lm([\alpha''_i|\beta''_i]_{Z_\gamma\transpose Z_\gamma})
\]
by Lemma \ref{lem:lm of stand monom} \ref{item:lm muwu} and Lemma \ref{lem:lm of zgamma t}
\ref{item:lm of zgamma t}.
The rest is proved along the same line with (1).
\end{proof}

\section{Doset Hibi rings}
\label{sec:doshibi}
%
%


\mysubsection{}
Let $P$ be a finite partially ordered set (poset for short).
The length of a chain (totally ordered subset) $X$ of $P$
is $|X|-1$.
The rank of $P$, denoted by $\rank P$, is the maximum
of the lengths of chains in $P$.
A poset is said to be
pure if its all maximal chains have the same length.
For $x$, $y\in P$, $y$ covers $x$, 
denoted by $x\covered y$,
means $x<y$ and 
there is no $z\in P$ such that $x<z<y$.
For $x$, $y\in P$ with $x\leq y$, we set
$[x,y]_P\define\{z\in P\mid x\leq z\leq y\}$.
The height
of an element $x\in P$,
denoted by $\height_Px$ or simply $\height x$ 
is the rank of $\{y\in P\mid y\leq x\}$.

\mysubsection{}
Let $H$ be a finite distributive lattice.
A \joinirred\ element 
in $H$
is an element $\alpha \in H$ such that $\alpha$ can not be expressed as a join of
two elements different from $\alpha$,
i.e.,  $\alpha=\beta\join\gamma\Rightarrow \alpha=\beta$ or $\alpha=\gamma$.
Note that we treat the unique minimal element of $H$ as a \joinirred\  element.

Let $P$ be the set of all \joinirred\ elements in $H$.
Then 
it is known
that $H$ is isomorphic to 
$J(P)\setminus\{\emptyset\}$ ordered by inclusion,
where $J(P)$ is  the set of all poset ideals of $P$.
The isomorphisms $\Phi\colon H\to J(P)\setminus\{\emptyset\}$ and
$\Psi\colon J(P)\setminus\{\emptyset\}\to H$ are given by
\[
\begin{array}{ll}
\Phi(\alpha)\define
	\{x\in P\mid x\leq\alpha\text{ in $H$}\}&\text{for $\alpha\in H$ and}
\\
\Psi(I)\define
	\displaystyle\bigvee_{x\in I}x&\text{for $I\in J(P)\setminus\{\emptyset\}$.}
\end{array}
\]


Let $P$ and $Q$ be posets.
A map $\nu \colon P\to Q$ is an order reversing map
or an anti-homomorphism
if $x\leq y$ in $P$ implies $\nu(x)\geq \nu(y)$ in $Q$
and strictly order reversing if
$x<y$ in $P$ implies $\nu(x)>\nu(y)$ in $Q$.
The set of all order reversing maps 
from $P$ to $\NNN$
is denoted by $\overline\TTTTT(P)$.
The set of all 
strictly order reversing maps
from $P$ to $\NNN\setminus\{0\}$
is denoted by 
$\TTTTT(P)$.

Let $H$ be a finite distributive lattice, $P$ the set of \joinirred\ elements
and $\{T_x\}_{x\in P}$ a family of indeterminates indexed by $P$.
We set $T^\nu=\prod_{x\in P}T_x^{\nu(x)}$ for $\nu\in \overline\TTTTT(P)$.
Hibi \cite{hib} defined the ring $\RRRRR_K(H)$, called the Hibi ring nowadays,
as follows.
\begin{definition}\rm
$\RRRRR_K(H)\define K[\prod_{x\leq \alpha}T_x\mid\alpha\in H]$.
\end{definition}
Hibi \cite{hib} showed the following theorems.
\begin{thm}
$\RRRRR_K(H)$ is an ASL over $K$ on $H$ with structure map
$\alpha\mapsto\prod_{x\leq \alpha}T_x$.
The straightening law of $\RRRRR_K(H)$ is 
$\alpha\beta=(\alpha\meet\beta)(\alpha\join\beta)$ for $\alpha$, $\beta\in H$
with $\alpha\not\sim\beta$.
\end{thm}
\begin{thm}
\label{thm:hibi basis}
$\RRRRR_K(H)$ is a $K$-vector space with basis
$\{T^\nu\mid\nu\in\overline\TTTTT(P)\}$.
In particular, by Theorem \ref{thm:hoc normal},
$\RRRRR_K(H)$ is a normal affine semigroupring
and hence \cm.
Further, by Theorem \ref{thm:sta can}, the canonical module of 
$\RRRRR_K(H)$ is a $K$-vector space with basis
$\{T^\nu\mid\nu\in\TTTTT(P)\}$.
\end{thm}

By analyzing the relations of the
leading monomials of 
the elements of $\Delta(m\times n;\gamma^*)$,
embeded by the structure map of the
ASL
$K[Z_\gamma]$ of Lemma \ref{lem:lm of stand monom} \ref{item:xwu isom},
we see the following
\begin{prop}
$\ini K[Z_\gamma]$ is the Hibi ring over $K$ on $\Delta(m\times n;\gamma^*)$.
In particular, it is normal.
In particular, by \cite[Corollary 2.3]{chv}, 
$K[Z_\gamma]$ 
is normal and \cm\ and
has rational singularities if $\chara K=0$
and is
F-rational if $\chara K>0$.
\end{prop}

\mysubsection{}
Let $L$ be another distributive lattice,
$Q$ the set of \joinirred\ elements in $L$ 
and
$\varphi\colon H\to L$ a surjective lattice homomorphism,
i.e.,
a surjective map from $H$ to $L$ such that
$\varphi(x\join y)=\varphi(x)\join\varphi(y)$
and
$\varphi(x\meet y)=\varphi(x)\meet\varphi(y)$.
We set
\[
\varphi^\ast(\beta)\define\bigmeet_{\alpha\in\varphi^{-1}(\beta)}\alpha.
\]
Then the following facts hold.
\begin{lemma}\label{lem:varphi ast basic}
\begin{enumerate}
\item\label{varphi ast basic 1}
$\varphi^\ast(\varphi(\alpha))\leq\alpha$ for any $\alpha\in H$.
\item\label{varphi ast basic 2}
$\varphi(\varphi^*(\beta))=\beta$ for any $\beta\in L$.
\item\label{varphi ast basic 3}
$\beta_1\leq \beta_2\Rightarrow \varphi^*(\beta_1)\leq\varphi^*(\beta_2)$.
\item\label{varphi ast basic 4}
For $\alpha\in H$ and $\beta\in L$,
$\varphi(\alpha)\geq\beta\iff\alpha\geq\varphi^*(\beta)$.
\item\label{varphi ast basic 6}
$y\in Q\Rightarrow \varphi^*(y)\in P$.
\item\label{varphi ast basic 7}
If $y_0$ is the unique minimal element of $L$,
then $\varphi^*(y_0)$ is the unique minimal element of $H$.
\end{enumerate}
\end{lemma}
\begin{proof}
\leavevmode
\ref{varphi ast basic 1}:
Since $\alpha\in \varphi^{-1}(\varphi(\alpha))$,
we see that
\[
\alpha\geq\bigmeet_{\alpha'\in\varphi^{-1}(\varphi(\alpha))}\alpha'=\varphi^*(\varphi(\alpha)).
\]
\ref{varphi ast basic 2}:
Since $\varphi$ is a lattice homomorphism, it holds
\[
\varphi(\varphi^\ast(\beta))
=\varphi(\bigmeet_{\alpha\in\varphi^{-1}(\beta)}\alpha)
=\bigmeet_{\alpha\in\varphi^{-1}(\beta)}\varphi(\alpha)
=\beta.
\]
\ref{varphi ast basic 3}:
Take $\alpha_0\in\varphi^{-1}(\beta_1)$.
Then for any $\alpha\in\varphi^{-1}(\beta_2)$,
$\varphi(\alpha_0\meet\alpha)=\varphi(\alpha_0)\meet\varphi(\alpha)=\beta_1\meet\beta_2=\beta_1$.
Therefore, $\alpha_0\meet\alpha\in\varphi^{-1}(\beta_1)$
so $\varphi^\ast(\beta_1)=\meet_{\alpha'\in\varphi^{-1}(\beta_1)}\alpha'
\leq \alpha_0\meet\alpha\leq \alpha$.
Since $\alpha$ is an arbitrary element of $\varphi^{-1}(\beta_2)$,
we see that
\[
\varphi^\ast(\beta_1)\leq\bigmeet_{\alpha\in\varphi^{-1}(\beta_2)}\alpha=\varphi^\ast(\beta_2).
\]
\ref{varphi ast basic 4}:
If $\varphi(\alpha)\geq \beta$, then $\alpha\geq\varphi^*(\varphi(\alpha))\geq\varphi^*(\beta)$
by \ref{varphi ast basic 1} and \ref{varphi ast basic 3}.
On the other hand, if $\alpha\geq\varphi^*(\beta)$, then
$\varphi(\alpha)\geq\varphi(\varphi^*(\beta))=\beta$ by \ref{varphi ast basic 2}.

\noindent\ref{varphi ast basic 6}:
If $\varphi^*(y)=\alpha\join\alpha'$, then
$\varphi(\alpha)\join\varphi(\alpha')=\varphi(\alpha\join\alpha')=\varphi(\varphi^\ast(y))=y$.
Therefore  $y=\varphi(\alpha)$ or $y=\varphi(\alpha')$,
since $y$ is \joinirred.
By symmetry, we may assume that $y=\varphi(\alpha)$.
Then $\varphi^\ast(y)=\varphi^\ast(\varphi(\alpha))\leq\alpha$ by \ref{varphi ast basic 1}.
On the other hand,
$\varphi^\ast(y)=\alpha\join\alpha'\geq \alpha$.
Therefoer $\varphi^\ast(y)=\alpha$.

\noindent\ref{varphi ast basic 7}:
Clear.
\end{proof}

By \ref{varphi ast basic 2} and \ref{varphi ast basic 3}  of Lemma \ref{lem:varphi ast basic}, 
we may regard $L$ as a subposet 
(not a sublattice) of $H$ by $\varphi^*$.
Then, by \ref{varphi ast basic 6}  of Lemma \ref{lem:varphi ast basic}, 
$Q$ is identified to a subset of $P$.
Note that $Q$ contains the unique minimal element of $H$ under this identification
by \ref{varphi ast basic 7} of Lemma \ref{lem:varphi ast basic}.
\begin{lemma}\label{lem:varphi iden}
By the identification above, the composition map
$J(P)\setminus\{\emptyset\}\simeq H\stackrel{\varphi}{\to} L\simeq J(Q)\setminus\{\emptyset\}$
is identical to the map $I\mapsto I\cap Q$.
\end{lemma}
\begin{proof}
Recall that lattice isomorphisms 
$H\simeq J(P)\setminus\{\emptyset\}$ and 
$L\simeq J(Q)\setminus\{\emptyset\}$
are obtained by
$\Phi_H\colon H\to J(P)\setminus\{\emptyset\}$, ($\alpha\mapsto\{x\in P\mid x\leq\alpha\}$) and 
$\Phi_L\colon L\to J(Q)\setminus\{\emptyset\}$, ($\beta\mapsto\{y\in Q\mid y\leq\beta\}$).
Let $\alpha$ be an arbitrary element of $H$.
Then $\Phi_L(\varphi(\alpha))=\{y\in Q\mid y\leq \varphi(\alpha)\}$.
Since $y\leq \varphi(\alpha)$ if and only if $\varphi^*(y)\leq \alpha$ by
\ref{varphi ast basic 4} of Lemma \ref{lem:varphi ast basic},
we see that by the embedding of $L$ in $H$ by $\varphi^*$,
$\Phi_L(\varphi(\alpha))$ is identified to
$\{y\in Q\mid y\leq \alpha\}$
which is equal to
$\{x\in P\mid x\leq \alpha\}\cap Q=\Phi_H(\alpha)\cap Q$.
\end{proof}

Let us recall the definition of dosets \cite[Section 18]{dep2}.
\begin{definition}\rm
A subset $D$ of $H\times H$ is called a doset of $H$ if
\begin{enumerate}
\item
$\{(\alpha,\alpha)\mid\alpha\in H\}\subset D\subset
\{(\alpha,\beta)\in H\times H\mid
\alpha\leq\beta\}$ and
\item
for $\alpha$, $\beta$ and $\gamma\in H$ with $\alpha\leq\beta\leq\gamma$,
\[
(\alpha,\gamma)\in D\iff(\alpha,\beta)\in D,(\beta,\gamma)\in D.
\]
\end{enumerate}
\end{definition}
Now let $L$ be another distributive lattice, 
$Q$ the set of \joinirred\ elements of $L$
and $\varphi\colon H\to L$ a
surjective lattice homomorphism.
We set
$D\define\{(\alpha,\beta)\in H\times H\mid
\alpha\leq\beta$, $\varphi(\alpha)=\varphi(\beta)\}$.
Then it is easily verified that $D$ is a doset of $H$.

We define the doset Hibi ring  as follows.
\begin{definition}\rm
The doset Hibi ring over $K$ defined by $\varphi$,
denoted by $\DDDDD_K(\varphi)$ is the $K$-subalgebra of $\RRRRR_K(H)$ generated
by $\{\alpha\beta\mid(\alpha,\beta)\in D\}$.
\end{definition}
Note that if $(\alpha,\beta)$, $(\alpha',\beta')\in D$, then
\begin{equation}
\begin{array}{rl}
&(\alpha\beta)({\alpha'}{\beta'})\\
=&({\alpha\meet\alpha'})({\alpha\join\alpha'})
	({\beta\meet\beta'})({\beta\join\beta'})\\
=& ({\alpha\meet\alpha'})
	({(\alpha\join\alpha')\meet(\beta\meet\beta')})
	({(\alpha\join\alpha')\join(\beta\meet\beta')})
	({\beta\join\beta'})
\end{array}
\label{eqn:doshibi str}
\end{equation}
and
\begin{eqnarray*}
\varphi((\alpha\join\alpha')\meet(\beta\meet\beta'))
&=&\varphi(\alpha\meet\alpha')\\
\varphi((\alpha\join\alpha')\join(\beta\meet\beta'))&=&
\varphi(\beta\join\beta'),
\end{eqnarray*}
since $\varphi(\alpha)=\varphi(\beta)$ and $\varphi(\alpha')=\varphi(\beta')$.
Therefore
$(\alpha\meet\alpha',(\alpha\join\alpha')\meet(\beta\meet\beta'))$,
$((\alpha\join\alpha')\join(\beta\meet\beta'),\beta\join\beta')\in D$.

In particular,
by applying the straightening law repeatedly, we see that
$\{
{\alpha_1}{\alpha_2}\cdots {\alpha_{2r-1}}{\alpha_{2r}}\mid
\alpha_1\leq\alpha_2\leq\cdots\leq\alpha_{2r}$,
$(\alpha_{2i-1},\alpha_{2i})\in D$ for $i=1$, \ldots, $r\}$
is a $K$-free basis of the doset Hibi ring.

\begin{example}\rm
\label{ex:extreme}
\begin{enumerate}
\item
\label{item:vero}
If $L$ is a set with only one element and $\varphi\colon H\to L$ is the unique map,
then $\DDDDD_K(\varphi)$ is the second Veronese subring of $\RRRRR_K(H)$.
\item
If $L=H$ and $\varphi$ is the identity map, then 
$\DDDDD_K(\varphi)=K[\alpha^2\mid\alpha\in H]$,
which is isomorphic to $\RRRRR_K(H)$.
\end{enumerate}
\end{example}

Regarding $\RRRRR_K(H)$ as a subring of $K[T_x\mid x\in P]$,
a standard monomial $\alpha_1\alpha_2\cdots\alpha_s$ with $\alpha_1\leq\alpha_2\leq\cdots\leq\alpha_s$
corresponds to $T^\nu$, where 
$\nu(x)=|\{i\mid x\leq\alpha_i\}|$.
On the other hand,
by the identification $L\simeq J(Q)\setminus\{\emptyset\}$ and $Q\subset P$,
$\varphi(\alpha)$ corresponds to 
$\{y\in Q\mid y\leq \alpha$ in $H\}$
by Lemma \ref{lem:varphi iden}.
Therefore a standard monomial $\alpha_1\alpha_2\cdots\alpha_{2r}$ on $H$
with $\alpha_1\leq\alpha_2\leq\cdots\leq\alpha_{2r}$ satisfies
$(\alpha_{2i-1},\alpha_{2i})\in D$ for $i=1$, \ldots, $r$ if and only if
$\nu(y)\equiv 0\pmod 2$ for any $y\in Q$,
where $\nu$ is an element of 
$\overline\TTTTT(P)$ corresponding to $\alpha_1\cdots\alpha_{2r}$.
Therefore we have the following
%
%
\begin{thm}\label{thm:dkvarphi descri}
$\DDDDD_K(\varphi)$ is a free $K$-module with basis
$\{T^\nu\mid\nu\in\overline\TTTTT(P)$,
$\nu(y)\equiv0\pmod 2$ for any $y\in Q\}$.
In particular,
by Theorem \ref{thm:hoc normal},
$\DDDDD_K(\varphi)$ is a normal affine semigroupring,
and therefore, is \cm.
\end{thm} 

Note the description of $\DDDDD_K(\varphi)$ in the theorem above depends  only
on $H$ and $Q$.
Therefore we can generalize the notion of doset Hibi ring as follows.
\begin{definition}\rm
Let $H$ be a finite distributive lattice,
$P$ the set of \joinirred\ elements of $H$,
$Q$ a subset of $P$.
We define the generalized doset Hibi ring defined by $H$ and $Q$,
denoted by $\DDDDD_K(H,Q)$,
as
$\DDDDD_K(H,Q)\define K[T^\nu\mid\nu\in\overline\TTTTT(P)$,
$\nu(y)\equiv0\pmod2$ for any $y\in Q]$.
\end{definition}
The following theorem is a direct consequence of the definition
and Theorem~\ref{thm:hoc normal}.
\begin{thm}\label{thm:gdoshibi is normal}
The generalized doset Hibi ring 
$\DDDDD_K(H,Q)$ 
is a free $K$-module with basis
$\{T^\nu\mid\nu\in\overline\TTTTT(P)$, $\nu(y)\equiv0\pmod2$ for any $y\in Q\}$.
In particular, $\DDDDD_K(H,Q)$ is a normal affine semigroupring
and hence is \cm.
\end{thm}
Note that if $Q$ contains the minimal element $x_0$ of $H$,
then by setting $L=J(Q)\setminus\{\emptyset\}$ and $\varphi$ the map
corresponding to $J(P)\setminus\{\emptyset\}\to J(Q)\setminus\{\emptyset\}$,
($I\mapsto I\cap Q$),
$\DDDDD_K(H,Q)$ is equal to $\DDDDD_K(\varphi)$.
In particular, $\DDDDD_K(H,Q)$ is a doset Hibi ring.

Here we consider a generating system of $\DDDDD_K(H,Q)$ as a  $K$-subalgebra of $\RRRRR_K(H)$.
Set $H_1\define\{\alpha\in H\mid y\leq\alpha$ for some $y\in Q\}$, $H_2=H\setminus H_1$
and 
$\psi\colon H_1\to J(Q)\setminus\{\emptyset\}$, ($\alpha\mapsto\{y\in Q\mid y\leq \alpha\}$).
\begin{lemma}\label{lem:in dkhq cond}
Let $\alpha_1\alpha_2\cdots\alpha_r$ be a standard monomial on $H$
with $\alpha_1\leq\alpha_2\leq\cdots\leq\alpha_r$,
$\alpha_1$, \ldots, $\alpha_s\in H_2$ and 
$\alpha_{s+1}$, \ldots, $\alpha_r\in H_1$.
Then $\alpha_1\cdots\alpha_r\in \DDDDD_K(H,Q)$ if and only if
$r-s\equiv0\pmod2$ and
$\psi(\alpha_{s+2i-1})=\psi(\alpha_{s+2i})$ for $i=1$, $2$, \ldots, $(r-s)/2$.
\end{lemma}
\begin{proof}
Let $\nu$ be the element of $\overline\TTTTT(P)$ corresponding to
$\alpha_1\cdots\alpha_r$.
Then $\nu(x)=|\{j\mid x\leq\alpha_j\}|$.
Therefore $\nu(y)=|\{j\mid j>s$ and $\psi(\alpha_j)\ni y\}|$ for $y\in Q$.
Since $\emptyset\neq\psi(\alpha_{s+1})\subset\psi(\alpha_{s+2})\subset\cdots\subset\psi(\alpha_r)$,
$\nu(y)\equiv0\pmod2$ for any $y\in Q$ if and only if
$r-s\equiv0\pmod2$ and
$\psi(\alpha_{s+2i-1})=\psi(\alpha_{s+2i})$
for $i=1$, $2$, \ldots, $(r-s)/2$.
The result follows from Theorem \ref{thm:gdoshibi is normal}.
\end{proof}

Now assume that $Q$ has a unique minimal element $y_0$.
Then $H_1=\{\alpha\in H\mid \alpha\geq y_0\}$ and $H_2=\{\alpha\in H\mid \alpha\not\geq y_0\}$.
\begin{lemma}\label{lem:meet join h1 h2}
\begin{enumerate}
\item
If $\alpha$, $\beta\in H_1$, then $\alpha\meet\beta$, $\alpha\join\beta\in H_1$.
\item
If $\alpha\in H_1$  and $\beta\in H_2$, then $\alpha\meet\beta\in H_2$, $\alpha\join\beta\in H_1$
and $\psi(\alpha\join\beta)=\psi(\alpha)$.
\item
If $\alpha$, $\beta\in H_2$, then $\alpha\meet\beta$, $\alpha\join\beta\in H_2$.
\item
If $\alpha$, $\beta\in H_1$, then $\psi(\alpha\meet\beta)=\psi(\alpha)\cap\psi(\beta)$
and $\psi(\alpha\join\beta)=\psi(\alpha)\cup\psi(\beta)$.
\item
If $\alpha$, $\beta$, $\alpha'$, $\beta'\in H_1$,
$\psi(\alpha)=\psi(\beta)$ and
$\psi(\alpha')=\psi(\beta')$,
then 
$\psi(\alpha\meet\alpha')=\psi((\alpha\join\alpha')\meet(\beta\meet\beta'))$ and
$\psi(\beta\join\beta')=\psi((\alpha\join\alpha')\join(\beta\meet\beta'))$.
\end{enumerate}
\end{lemma}
\begin{proof}
Let $\delta$ be an element of $H$ and $I$ the corresponding element of $J(P)\setminus\{\emptyset\}$
by the lattice isomorphism $H\simeq J(P)\setminus\{\emptyset\}$.
Then $\delta\in H_1$ if and only if $I\ni y_0$,
and $\meet$ and $\join$ operations in $H$ correspond to $\cap$ and $\cup$ in $J(P)\setminus\{\emptyset\}$
respectively.
Furthermore, $\psi(\delta)=I\cap Q$ if $\delta\in H_1$.
The lemma is proved straightfowardly by these facts.
\end{proof}
Now we state the following
\begin{prop}\label{prop:gen of gdoshibi}
If $Q$ has a unique minimal element, then
$\DDDDD_K(H,Q)$ is a $K$-subalgebra of $\RRRRR_K(H)$ generated by
$\{\alpha\beta\mid\alpha, \beta\in H_1$, $\alpha\leq\beta$  and  $\psi(\alpha)=\psi(\beta)\}\cup H_2$.
\end{prop}
\begin{proof}
We set $G=\{\alpha\beta\mid\alpha, \beta\in H_1$, $\alpha\leq\beta$  and  $\psi(\alpha)=\psi(\beta)\}\cup H_2$.
Then by Lemma \ref{lem:in dkhq cond}, we see that $\DDDDD_K(H,Q)\subset K[G]$.

To prove the reverse inclusion, let $\mu$ be an arbitrary product of elements of $G$.
Then by Lemma \ref{lem:meet join h1 h2} and equation \refeq{eqn:doshibi str}, we see that,
by repeated application of the straightening law of $\RRRRR_K(H)$,
$\mu$ is equal to a standard monomial $\alpha_1\cdots\alpha_r$ with
$\alpha_1\leq\cdots\leq\alpha_r$ such that for some $s$,
$\alpha_1$, \ldots, $\alpha_s\in H_2$,
$\alpha_{s+1}$, \ldots, $\alpha_r\in H_1$,
$r-s\equiv0\pmod2$ and
$\psi(\alpha_{s+2i-1})=\psi(\alpha_{s+2i})$ for $i=1$, $2$, \ldots, $(r-s)/2$.
Then again by Lemma \ref{lem:in dkhq cond},
we see that $\alpha_1\cdots\alpha_r\in\DDDDD_K(H,Q)$.
Therefore $\mu\in\DDDDD_K(H,Q)$.
\end{proof}

Now we apply the theory of (generalized) doset Hibi rings to the ring of 
absolute $\orth(m)$
and $\sorth(m)$-invariants.

Set $L=\{1,2,\ldots, m\}$ with reverse order and
$\varphi\colon\Gamma'(m\times n;\gamma)\to L$, 
($\delta\mapsto\msize\delta$).
Then $\varphi$ is a surjective lattice homomorphism
and $\varphi^*(j)=[b_1,b_2,\ldots,b_j]$ for $1\leq j\leq m$.
By analyzing the relations of the leading monomials of the sagbi basis
$\{[\alpha|\beta]_{Z_\gamma\transpose  Z_\gamma}\mid
(\alpha,\beta)\in D_\gamma\}$
of
$K[Z_\gamma]^{\orth(m,-)}=
K[Z_\gamma\transpose  Z_\gamma]$,
and the sagbi basis 
$\{[\alpha|\beta]_{Z_\gamma\transpose  Z_\gamma}\mid
(\alpha,\beta)\in D_\gamma$, $\msize\alpha<m\}
\cup\{\delta_{Z_\gamma}\mid
\delta\in\Gamma(m\times n;\gamma)\}$
of
$K[Z_\gamma]^{\sorth(m,-)}=
K[Z_\gamma\transpose  Z_\gamma, \Gamma(Z_\gamma)]$
(c.f.\ Theorem \ref{thm:sagbi of inv}),
we see the following
\begin{thm}\label{thm:ini is dhibi}
\begin{enumerate}
\item
$\ini K[Z_\gamma]^{\orth(m,-)}=
\ini K[Z_\gamma\transpose  Z_\gamma]$
is the doset Hibi ring  over $K$ defined by $\varphi$.
\item
$\ini K[Z_\gamma]^{\sorth(m,-)}=
\ini K[Z_\gamma\transpose  Z_\gamma, \Gamma(Z_\gamma)]$
is the generalized doset Hibi ring defined by
$\Gamma'(m\times n;\gamma)$
and $\{[b_1,\ldots,b_{m-1}], \ldots, [b_1,b_2], [b_1]\}$.
\end{enumerate}
In particular,
$\ini K[Z_\gamma]^{\orth(m,-)}$ and $\ini K[Z_\gamma]^{\sorth(m,-)}$ are 
normal affine semigrouprings and therefore
$K[Z_\gamma]^{\orth(m,-)}$ and $K[Z_\gamma]^{\sorth(m,-)}$ are normal and \cm, and
has rational singularities if $\chara K=0$
and is
F-rational if $\chara K>0$.
\end{thm}

\section{\gor\ property}
\label{sec:gor}

In this section, we stat a criterion of the \gor\  property of a
(generalized) doset Hibi ring.

Let $H$ be a finite distributive lattice, $P$ the set of \joinirred\ elements of 
$H$ and $Q$ a subset of $P$
and $\{T_x\}_{x\in P}$ the family of indeterminates indexed by $P$.
Set $\overline\TTTTT(P,Q)\define\{\nu\in\overline\TTTTT(P)\mid \nu(y)\equiv0\pmod2$ for any $y\in Q\}$
and $\TTTTT(P,Q)\define\{\nu\in\TTTTT(P)\mid \nu(y)\equiv0\pmod2$ for any $y\in Q\}$.
Then since $\DDDDD_K(H,Q)=\bigoplus_{\nu\in\overline\TTTTT(P,Q)} K T^\nu$,
we see, by Theorem \ref{thm:sta can}, that the canonical module of
$\DDDDD_K(H,Q)$ is $\bigoplus_{\nu\in\TTTTT(P,Q)} K T^\nu$.

We define a order on $\TTTTT(P,Q)$ by
$\nu\leq\nu'\stackrel{\rm def}{\iff}\nu'-\nu\in\overline\TTTTT(P,Q)$,
where $(\nu'-\nu)(x)\define\nu'(x)-\nu(x)$ for any $x\in P$.
Then it is clear that $T^\nu$ is a generator of the canonical module of 
$\DDDDD_K(H,Q)$ if and only if $\nu$ is a minimal element of $\TTTTT(P,Q)$.
In particular, $\DDDDD_K(H,Q)$ is \gor\ if and only if $\TTTTT(P,Q)$ has a
unique minimal element.

Now we set $P^+\define P\cup\{\infty\}$, where $\infty$ is a new element
and $x<\infty$ for any $x\in P$ and $Q^+\define Q\cup\{\infty\}$.
We also set
$\tilde P\define P^+\cup\{(y_1,y_2)\in Q^+\times Q^+\mid y_1\covered y_2$ in $P^+\}$
and define the order on $\tilde P$ by extending the order relation on $P$ by
$y_1<(y_1,y_2)<y_2$ for $y_1$, $y_2\in Q^+$ with $y_1\covered y_2$ in $P^+$.
Then we have the following result.
\begin{thm}\label{thm:gor cri}
$\DDDDD_K(H,Q)$ is \gor\ if and only if the following two conditions are satisfied.
\begin{enumerate}
\item\label{item:tilde p pure}
$\tilde P$ is pure.
\item\label{item:rank even}
For any $y$, $y'\in Q^+$ with $y<y'$,
$\rank[y,y']_{\tilde P}\equiv0\pmod2$.
\end{enumerate}
\end{thm}
\begin{proof}
First we define $\tilde\nu_0\colon\tilde P\to\NNN$ by
downward induction as follows.
\begin{itemize}
\item
$\tilde\nu_0(\infty)=0$.
\item
If $x\neq \infty$ and $\tilde\nu_0(x')$ is defined for any $x'\in\tilde P$ with $x'>x$,
we set
\[
\tilde\nu_0(x)=
\left\{
\begin{array}{ll}
\max\{\tilde\nu_0(x')+1\mid x'>x\}&\mbox{ if $x\not\in Q$,}
\\
2\lceil\max\{\tilde\nu_0(x')+1\mid x'>x\}/2\rceil&\mbox{ if $x\in Q$.}
\end{array}
\right.
\]
\end{itemize}
Then it is clear that 
$\nu_0\define\tilde\nu_0|_P\in\TTTTT(P,Q)$ and
$\nu(x)\geq\nu_0(x)$ for any $\nu\in\TTTTT(P,Q)$ and $x\in P$.
In particular, $\nu_0$
is a minimal element of $\TTTTT(P,Q)$.

Now we start the proof of the theorem and 
we first prove the ``if'' part of the theorem.
By assumptions (1) and (2), it is easily verified that
$\tilde\nu_0(x)=\tilde\nu_0(x')+1$ for any $x$, $x'\in \tilde P$ with $x\covered x'$.
Therefore, by the definition of $\tilde P$, we see that
$\nu-\nu_0\in\overline\TTTTT(P,Q)$ for any $\nu\in\TTTTT(P,Q)$.
Thus, $\nu_0$ is the unique minimal element of $\TTTTT(P,Q)$ and $\DDDDD_K(H,Q)$ is \gor.

Next we prove the ``only if'' part.
If (1) or (2) is not valid, then there are $x$, $x'\in \tilde P$ with
$x\covered x'$ in $\tilde P$ such that $\tilde\nu_0(x)\geq\tilde\nu_0(x')+2$.
Since $x'\neq\infty$ in such a case,  we see that there are $x$, $x'\in P$ 
such that $x\covered x'$ in $P$,
$\{x,x'\}\not\subset Q$ and
$\nu_0(x)\geq\nu_0(x')+2$
or
there are $y$, $y'\in Q$ such that $y\covered y'$ in $P$ and $\nu_0(y)\geq\nu_0(y')+4$.
We classify this phenomena into three cases.
\begin{case}\label{case:gor cri 1}\rm
There are $y$, $y'\in Q$ such that $y\covered y'$ in $P$ and 
$\nu_0(y)\geq\nu_0(y')+4$.
\end{case}
\begin{case}\label{case:gor cri 2}\rm
There are $x$, $x'\in P$ such that $x\covered x'$, $x\not\in Q$ and
$\nu_0(x)\geq\nu_0(x')+2$.
\end{case}
\begin{case}\label{case:gor cri 3}\rm
There are $x$, $x'\in P$ such that $x\covered x'$, $x\in Q$, $x'\not\in Q$ and
$\nu_0(x)\geq\nu_0(x')+2$.
\end{case}

As for Case \ref{case:gor cri 1}, we set
\[
\nu(z)=
\left\{
\begin{array}{ll}
\nu_0(z)&\mbox{ $z\not\leq y'$ or $z=y$,}
\\
\nu_0(z)+2&\mbox{ $z\leq y'$ and $z\neq y$.}
\end{array}
\right.
\]
Then it is easily verified that $\nu\in\TTTTT(P,Q)$ and $\nu-\nu_0\not\in\overline\TTTTT(P,Q)$
since $(\nu-\nu_0)(y)=0$ and $(\nu-\nu_0)(y')=2$.
Therefore, there is a minimal element of $\TTTTT(P,Q)$ other than $\nu_0$.

As for Case \ref{case:gor cri 2}, we set
\[
\nu(z)=
\left\{
\begin{array}{ll}
\nu_0(z)&\mbox{ $z\not\leq x'$,}\\
\nu_0(z)+1&\mbox{ $z=x$,}\\
\nu_0(z)+2&\mbox{ $z\leq x'$ and $z\neq x$.}
\end{array}
\right.
\]
Then it is easily verified that $\nu\in\TTTTT(P,Q)$ and $\nu-\nu_0\not\in\overline\TTTTT(P,Q)$
since $(\nu-\nu_0)(x)=1$ and $(\nu-\nu_0)(x')=2$.
Therefore, there is a minimal element of $\TTTTT(P,Q)$ other than $\nu_0$.

Finally,
as for Case \ref{case:gor cri 3}, we set
\[
\nu(z)=
\left\{
\begin{array}{ll}
\nu_0(z)&\mbox{ $z\not\leq x'$ or $z=x$,}\\
\nu_0(z)+1&\mbox{ $z=x'$,}\\
\nu_0(z)+2&\mbox{ $z< x'$ and $z\neq x$.}
\end{array}
\right.
\]
Again it is easily verified that $\nu\in\TTTTT(P,Q)$ and $\nu-\nu_0\not\in\overline\TTTTT(P,Q)$.
Therefore, there is a minimal element of $\TTTTT(P,Q)$ other than $\nu_0$.

In any case, we see that $\DDDDD_K(H,Q)$ is not \gor.
\end{proof}

\begin{example}\rm
\begin{enumerate}
\item
If
\[
P=
\vcenter{\hsize=40\unitlength
\begin{picture}(40,50)
\put(20,10){\circle*{4}}
\put(10,20){\circle*{2}}
\put(10,30){\circle*{4}}
\put(10,40){\circle*{2}}
\put(30,30){\circle*{4}}

\put(10,20){\line(0,1){20}}
\put(10,20){\line(1,-1){10}}
\put(20,10){\line(1,2){10}}
\end{picture}},
\]
where big dots express elements  of $Q$, then
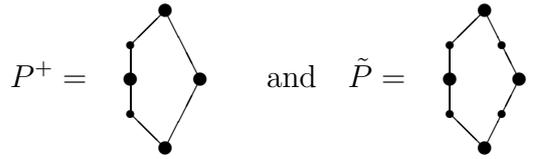
\[
P^+=
\vcenter{\hsize=40\unitlength
\begin{picture}(40,60)
\put(20,10){\circle*{4}}
\put(10,20){\circle*{2}}
\put(10,30){\circle*{4}}
\put(10,40){\circle*{2}}
\put(30,30){\circle*{4}}

\put(20,50){\circle*{4}}

\put(10,20){\line(0,1){20}}
\put(10,20){\line(1,-1){10}}
\put(20,10){\line(1,2){10}}

\put(10,40){\line(1,1){10}}
\put(20,50){\line(1,-2){10}}
\end{picture}}
\quad\mbox{and}\quad
\tilde P=
\vcenter{\hsize=40\unitlength
\begin{picture}(40,60)
\put(20,10){\circle*{4}}
\put(10,20){\circle*{2}}
\put(10,30){\circle*{4}}
\put(10,40){\circle*{2}}
\put(30,30){\circle*{4}}

\put(20,50){\circle*{4}}
\put(25,20){\circle*{2}}
\put(25,40){\circle*{2}}

\put(10,20){\line(0,1){20}}
\put(10,20){\line(1,-1){10}}
\put(20,10){\line(1,2){10}}

\put(10,40){\line(1,1){10}}
\put(20,50){\line(1,-2){10}}
\end{picture}}.
\]
Therefore, $\DDDDD_K(H,Q)$ is \gor.
\item
If
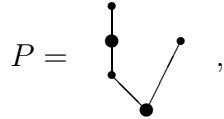
\[
P=
\vcenter{\hsize=40\unitlength
\begin{picture}(40,50)
\put(20,10){\circle*{4}}
\put(10,20){\circle*{2}}
\put(10,30){\circle*{4}}
\put(10,40){\circle*{2}}
\put(30,30){\circle*{2}}

\put(10,20){\line(0,1){20}}
\put(10,20){\line(1,-1){10}}
\put(20,10){\line(1,2){10}}
\end{picture}},
\]
then
\[
P^+=
\vcenter{\hsize=40\unitlength
\begin{picture}(40,60)
\put(20,10){\circle*{4}}
\put(10,20){\circle*{2}}
\put(10,30){\circle*{4}}
\put(10,40){\circle*{2}}
\put(30,30){\circle*{2}}

\put(20,50){\circle*{4}}

\put(10,20){\line(0,1){20}}
\put(10,20){\line(1,-1){10}}
\put(20,10){\line(1,2){10}}

\put(10,40){\line(1,1){10}}
\put(20,50){\line(1,-2){10}}
\end{picture}}
=\tilde P.
\]
Therefore, $\tilde P$ is not pure and $\DDDDD_K(H,Q)$ is not \gor.
\item
If
\[
P=
\vcenter{\hsize=40\unitlength
\begin{picture}(40,70)
\put(20,10){\circle*{4}}
\put(10,20){\circle*{2}}
\put(10,30){\circle*{2}}
\put(10,40){\circle*{4}}
\put(10,50){\circle*{2}}
\put(10,60){\circle*{2}}
\put(30,30){\circle*{4}}
\put(30,50){\circle*{4}}

\put(10,20){\line(0,1){40}}
\put(10,20){\line(1,-1){10}}
\put(20,10){\line(1,2){10}}
\put(30,30){\line(0,1){20}}
\end{picture}},
\]
then
\[
P^+=
\vcenter{\hsize=40\unitlength
\begin{picture}(40,70)
\put(20,10){\circle*{4}}
\put(10,20){\circle*{2}}
\put(10,30){\circle*{2}}
\put(10,40){\circle*{4}}
\put(10,50){\circle*{2}}
\put(10,60){\circle*{2}}
\put(30,30){\circle*{4}}
\put(30,50){\circle*{4}}

\put(10,20){\line(0,1){40}}
\put(10,20){\line(1,-1){10}}
\put(20,10){\line(1,2){10}}
\put(30,30){\line(0,1){20}}

\put(20,70){\circle*{4}}
\put(10,60){\line(1,1){10}}
\put(20,70){\line(1,-2){10}}
\end{picture}}
\quad\mbox{and}\quad
\tilde P=
\vcenter{\hsize=40\unitlength
\begin{picture}(40,70)
\put(20,10){\circle*{4}}
\put(10,20){\circle*{2}}
\put(10,30){\circle*{2}}
\put(10,40){\circle*{4}}
\put(10,50){\circle*{2}}
\put(10,60){\circle*{2}}
\put(30,30){\circle*{4}}
\put(30,50){\circle*{4}}

\put(10,20){\line(0,1){40}}
\put(10,20){\line(1,-1){10}}
\put(20,10){\line(1,2){10}}
\put(30,30){\line(0,1){20}}

\put(20,70){\circle*{4}}
\put(10,60){\line(1,1){10}}
\put(20,70){\line(1,-2){10}}

\put(25,20){\circle*{2}}
\put(25,60){\circle*{2}}
\put(30,40){\circle*{2}}

\end{picture}}.
\]
Thus, there are $y$, $y'\in Q^+$ such that 
$y<y'$ and $\rank[y,y']_{\tilde P}=3$.
Therefore, $\DDDDD_K(H,Q)$ is not \gor.
\end{enumerate}

\end{example}

Finally, we apply Theorem \ref{thm:gor cri} to obtain a criterion of \gor\
property of $\ini K[Z_\gamma]^{\orth(m,-)}$ and $\ini K[Z_\gamma]^{\sorth(m,-)}$.
Recalled that we have fixed $\gamma=[b_1, \ldots, b_m]\in\Gamma(m\times n)$.
Note by Corollary \ref{cor:same hilb},
$K[Z_\gamma]^{\orth(m,-)}$ (resp.\ $K[Z_\gamma]^{\sorth(m,-)}$) is \gor\
if and only if so is $\ini K[Z_\gamma]^{\orth(m,-)}$ (resp.\ $\ini K[Z_\gamma]^{\sorth(m,-)}$).

First we note the following lemma which is easily proved.
\begin{lemma}\label{lem:cri joinirred}
Let $H$ be a distributive lattice and $x$ an element of $H$ which is not minimal.
Then $x$ is \joinirred\ if and only if there uniquely exists an element $\alpha\in H$
such that $\alpha\covered x$.
\end{lemma}
In order to apply Theorem \ref{thm:gor cri} 
to the ring of absolute invariants, we first analyze the structure of the
set  of \joinirred\ elements of $\Gamma'(m\times n;\gamma)$.
Let $P$ be the set of \joinirred\ elements of $\Gamma'(m\times n;\gamma)$ and
$\delta=[d_1,d_2,\ldots,d_t]$ be an element of $\Gamma'(m\times n;\gamma)\setminus\{\gamma\}$.

We consider the condition that $\delta\in P$.
First we state the following lemma which is easily proved.
\begin{lemma}
\label{lem:cover}
Let $[c_1,\ldots, c_{s}]$,
$[c'_1,\ldots, c'_{s'}]\in\Gamma'(m\times n;\gamma)$.
Then 
$$
[c_1,\ldots, c_{s}]\covered[c'_1, \ldots, c'_{s'}]
$$
in $\Gamma'(m\times n;\gamma)$ if and only if 
\begin{enumerate}
\item
$s=s'+1$, $c_{s}=m$ and $c'_i=c_i$ for $1\leq i\leq s'$
or
\item
$s=s'$ and there exists $i$ such that $c'_i=c_i+1$ and
$c'_j=c_j$ for $j\neq i$.
\end{enumerate}
\end{lemma}
Now consider the condition of $\delta=[d_1,\ldots,d_t]\in P$.
First consider the case where $\msize \delta=m$.
By Lemmas \ref{lem:cri joinirred} and \ref{lem:cover},
we see that $\delta$ is \joinirred\ if and only if there
is a unique index $i$ such that 
$$
d_i>b_i\quad\mbox{and}\quad
d_i>d_{i-1}+1,
$$
where we set $d_0=0$.
Note that $i=\min\{j\mid d_j>b_j\}$
and $\delta=[b_1, \ldots, b_{i-1}, d_i, \max\{d_i+1,b_{i+1}\},\ldots,\max\{d_i+m-i,b_m\}]$
in this case.

Next consider the case where $\msize \delta=t<m$.
By Lemmas \ref{lem:cri joinirred} and \ref{lem:cover},
we see that $\delta$ is \joinirred\ if and only if $\delta=[b_1, \ldots, b_t]$ or
there is $i$ such that
$$
d_j=b_j\mbox{ for $j<i$ and }
d_j=n-t+j\mbox{  for $j\geq i$.}
$$

Now we construct a map from $P\setminus\{\gamma\}$ to $\ZZZ\times \NNN$.
Let $\delta=[d_1, \ldots, d_t]$
be an arbitrary element of $P\setminus\{\gamma\}$.
Set $i=\min\{j\mid d_j>b_j\}$, where we set $d_{t+1}=n+1$ in the case
where $\delta=[b_1, \ldots, b_t]$.
Set $\xi(\delta)=(n-d_i-m+i,i-1)$.
Then $\xi$ is a map from $P\setminus\{\gamma\}$ to $\ZZZ\times\NNN$.
Note that  $i$ is the unique index such that 
$$
d_i>b_i\quad\mbox{and}\quad
d_i>d_{i-1}+1,
$$
where we set $d_0=0$.
Also note that $n-d_i-m+i=|\{j\mid j\not\in\{d_1, \ldots, d_m\}$ and $d_i<j\leq n\}|$
in the case where $t=m$ and $n-d_i-m+i=t-m$ in the case where $t<m$.
In particular, $\msize \delta=m$ if and only if the first entry of
$\xi(\delta)$ is non-negative.
Note also that $i-1=|\{j\mid j\in\{d_1, \ldots, d_t\}$ and $1\leq j< d_i\}|$.

Suppose $\delta'=[d'_1,\ldots, d'_{t'}]\in P\setminus\{\gamma\}$.
Set $\imath'=\min\{j\mid d'_j>b_j\}$, where we set $d'_{t'+1}=n+1$ as above.
First assume that $\xi(\delta')=\xi(\delta)$.
Then $\imath'=i$, $d'_i=d_i$ and 
$d'_j=\max\{d'_i+j-i, b_j\}=
\max\{d_i+j-i, b_j\}=d_j$ 
for $i<j\leq m$ in the case where $t=m$
or 
$t'-m=t-m$ and
$d'_j=n-t'+j=n-t+j=d_j$ for $i\leq j\leq t$ in the case where $t<m$.
Therefore, we see that $\delta'=\delta$, i.e., $\xi$ is injective.

Next assume that $\delta'>\delta$.
Then $\imath'\leq i$, $d'_i\geq d_i$ and $d'_i-d'_{\imath'}=i-\imath'$,
since $\max\{d'_{\imath'}+i-\imath',b_i\}=d'_i\geq d_i>b_i$,
where we set $d'_{t'+j}=n+j$ for $j>0$.
Therefore,
$n-d'_{\imath'}-m+\imath'=n-d'_i-m+i\leq n-d_i-m+i$
and $\imath'-1\leq i-1$.
Thus, we see that $\xi$ is an anti-homomorphism from $P\setminus\{\gamma\}$
to $\ZZZ\times \NNN$.

Set $\xi(\delta)=(a,b)$.
First suppose that $b>0$.
Then $i\geq 2$ and 
$\delta''=[d_1, \ldots, d_{i-2}, d_i-1,d_i,\ldots, d_t]$
is an element of $P\setminus\{\gamma\}$ such that
$\xi(\delta'')=(a,b-1)$,
where we set $d_{t+1}=n+1$ in the case where $\delta=[b_1, \ldots, b_t]$
as above.
Next suppose that $a-b>-m$ and $(a,b)\neq(-m+1,0)$.
Since $\xi([b_1, \ldots, b_t])=(t-m,t)$,
we see that $\delta\neq[b_1, \ldots b_t]$, i.e.,
$i\leq t$.
Set $\delta''=[d_1, \ldots, d_{i-1}, d_i+1, \max\{d_i+2, d_{i+1}\},
\ldots, \max\{d_i+m-i+1, d_m\}]$
in the case where $t=m$ and $n-d_i>m-i$
or
$\delta''=[d_1,\ldots, d_{i-1}, d_{i+1}, \ldots d_t]$
in the case where $t<m$ or $n-d_i=m-i$.
Note also $\delta''=[d_1, \ldots, d_{i-1}]=[b_1, \ldots, b_{i-1}]$
in the case where $i=t$ and $d_t=n$.
It is easily verified that $\delta''\in P\setminus\{\gamma\}$
and $\xi(\delta'')=(a-1,b)$.

Set $P'=\{(x,y)\in\ZZZ\times\NNN\mid y\geq 0$, $x-y\geq -m$, $(x,y)\neq (-m,0)\}$.
Then $P'$ is a subposet of $\ZZZ\times \NNN$ with unique minimal element
$(-m+1,0)$.
Further, we see that $P\setminus\{\gamma\}$  is anti-isomorphic to 
the poset ideal $\image\xi$ of $P'$ by the above argument.
Since $\xi$ is an anti-homomorphism, $\image\xi$ is the poset ideal
of $P'$ generated by
$\{\xi(\delta)\mid\delta$ is a minimal element of $P\setminus\{\gamma\}\}$.

By Lemma \ref{lem:cri joinirred}, we see that $\delta\in\Gamma'(m\times n;\gamma)$
is a minimal element of $P\setminus\{\gamma\}$ if and only if $\delta\covers\gamma$
in $\Gamma'(m\times n;\gamma)$.
Thus, in order to study the set of minimal elements of $P\setminus\{\gamma\}$,
we consider the set 
$\{\delta\in\Gamma'(m\times n;\gamma)\mid\delta\covers\gamma\}$.

We separate $\{b_1,b_2,\ldots, b_m\}$ 
into blocks and define the gaps between them.
Set
$U=
\{u\mid 1\leq u\leq m$, $b_u+1<b_{u+1}\}=\{u(1),\ldots, u(k)\}$ with
$u(1)<\cdots<u(k)$, where we set $b_{m+1}\define n+1$.
Then $U\neq \emptyset$ if and only if $\gamma\neq[n-m+1,\ldots, n]$.

First consider the case where $\gamma=[n-m+1,\ldots, n]$.
Then $P\setminus\{\gamma\}$ has the unique minimal element $[n-m+1,\ldots, n-1]
=[b_1,\ldots, b_{m-1}]$.
Since $\xi([b_1,\ldots, b_{m-1}])=(-1,m-1)$,
we see that $P\setminus\{\gamma\}$ is anti-isomorphic to
$\{(x,y)\in\ZZZ\times\NNN\mid 0\leq y\leq m-1$, $x\leq -1$, $x-y\geq -m$, $(x,y)\neq (-m,0)\}$.
In particular, $P$ is pure.

Next consider the case where $\gamma\neq[n-m+1,\ldots, n]$.
$U\neq\emptyset$ and $k\geq 1$ in this case.
Set
$\chi_0\define\{1,2,\ldots,b_1-1\}$,
$B_1\define\{b_1,b_2,\ldots, b_{u(1)}\}$,
$\chi_1\define\{b_{u(1)}+1,b_{u(1)}+2,\ldots, b_{u(1)+1}-1\}$,
$B_2\define\{b_{u(1)+1},b_{u(1)+2},\ldots, b_{u(2)}\}$,
$\chi_2\define\{b_{u(2)}+1,b_{u(2)}+2,\ldots, b_{u(2)+1}-1\}$,
$B_3\define\{b_{u(2)+1},b_{u(2)+2},\ldots, b_{u(3)}\}$,
\ldots,
$B_k\define\{b_{u(k-1)+1},b_{u(k-1)+2},\ldots, b_{u(k)}\}$,
$\chi_k\define\{b_{u(k)}+1,b_{u(k)}+2,\ldots, b_{u(k)+1}-1\}$
and
$B_{k+1}\define\{b_{u(k)+1},b_{u(k)+2},\ldots, b_{m}\}$.
Note that if $b_m<n$, then $u(k)=m$ and $B_{k+1}=\emptyset$.

By Lemma \ref{lem:cover} and the above remark, we see
that the minimal elements of $P\setminus\{\gamma\}$ are
$\gamma_1=[b_1$, \ldots, $b_{u(1)-1},b_{u(1)}+1,b_{u(1)+1}$, \ldots, $b_m]$,
$\gamma_2=[b_1$, \ldots, $b_{u(2)-1},b_{u(2)}+1,b_{u(2)+1}$, \ldots, $b_m]$,
\ldots,
$\gamma_{k-1}=[b_1$, \ldots, $b_{u(k-1)-1},b_{u(k-1)}+1,b_{u(k-1)+1}$, \ldots, $b_m]$
and
$\gamma_k=[b_1$, \ldots, $b_{m-1}$, $b_m+1]$
in the case where $b_m<n$ and
$\gamma_1$, \ldots, $\gamma_{k-1}$,
$\gamma_k=[b_1$, \ldots, $b_{u(k)-1},b_{u(k)}+1,b_{u(k)+1}$, \ldots, $b_m]$
and
$\gamma_{k+1}=[b_1,b_2,\ldots,b_{m-1}]$
in the case where $b_m=n$.

Note that
$$|B_j|=u(j)-u(j-1)$$ for $2\leq j\leq k$ and $$|B_{k+1}|=m-u(k)$$
in the case where $b_m=n$.
Further, 
$$
|\chi_{j-1}|+|B_j|=b_{u(j)}-b_{u(j-1)}
$$
for $2\leq j\leq k$ and 
$$
|\chi_k|+|B_{k+1}|=b_m-b_{u(k)}=n-b_{u(k)}
$$
in the case where $b_m=n$.
Thus, we see that 
$$
|\chi_{j-1}|=b_{u(j)}-b_{u(j-1)}-u(j)+u(j-1)
$$
for $2\leq j\leq k$ and
$$
|\chi_k|=n-b_{u(k)}-m+u(k)
$$
in the case where $b_m=n$.

Note also that
$$
\xi(\gamma_j)=(n-b_{u(j)}-1-m+u(j),u(j)-1)
$$
for $1\leq j\leq k$ and
$$
\xi(\gamma_{k+1})=(-1,m-1)
$$
in the case where $b_m=n$.
Therefore, the difference between the second entries of 
$\xi(\gamma_j)$ and $\xi(\gamma_{j-1})$ is $|B_j|$ and the
difference between the first entries of $\xi(\gamma_j)$
and $\xi(\gamma_{j-1})$ is $|\chi_{j-1}|$ 
for $2\leq j\leq k$
(resp.\ $2\leq j\leq k+1$)
if $b_m<n$ (resp.\ $b_m=n$).
Since $P'$ is a poset with unique minimal element $(-m+1,0)$,
we see that $\height_{P'}(a,b)=a+b+m-1$ for any $(a,b)\in P'$.
Therefore, we see that
$\height_{P'}\xi(\gamma_j)=\height_{P'}\xi(\gamma_{j-1})$
if and only if $|B_j|=|\chi_{j-1}|$
for $2\leq j\leq k$
(resp.\ $2\leq j\leq k+1$)
if $b_m<n$ (resp.\ $b_m=n$).

\begin{example}\rm
\label{ex:p setminus gamma}
\begin{enumerate}
\item
If $m=7$, $n=15$ and $\gamma=[2,6,7,8,10,13,14]$, then
$k=4$, $u(1)=1$, $u(2)=4$, $u(3)=5$, $u(4)=7$ and
$\chi_0=\{1\}$,
$B_1=\{2\}$,
$\chi_1=\{3,4,5\}$,
$B_2=\{6,7,8\}$,
$\chi_2=\{9\}$,
$B_3=\{10\}$,
$\chi_3=\{11,12\}$,
$B_4=\{13,14\}$,
$\chi_4=\{15\}$
and
$B_5=\emptyset$.
$\gamma_1=[3,6,7,8,10,13,14]$, 
$\gamma_2=[2,6,7,9,10,13,14]$, 
$\gamma_3=[2,6,7,8,11,13,14]$ and 
$\gamma_4=[2,6,7,8,10,13,15]$.
Further,
$\xi(\gamma_1)=(6,0)$,
$\xi(\gamma_2)=(3,3)$,
$\xi(\gamma_3)=(2,4)$ and
$\xi(\gamma_4)=(0,6)$.
Therefore, 
the Hasse diagram of $P\setminus\{\gamma\}$ is the following.
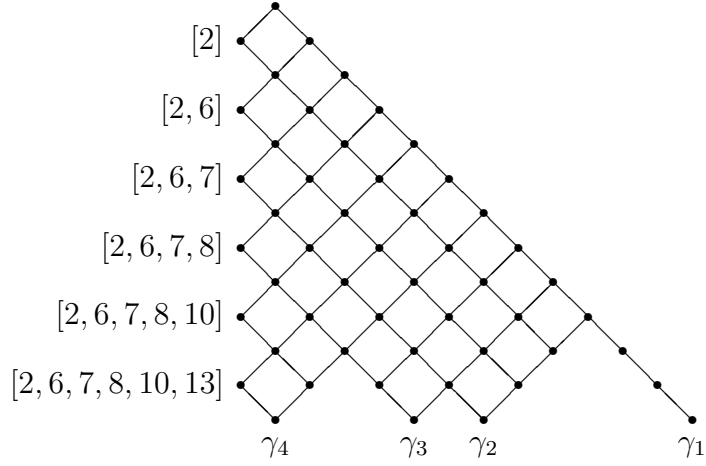
\begin{center}
\begin{picture}(160,140)(-20,0)

\multiput(20,130)(10,-10){13}{\circle*{2}} \put(20,130){\line(1,-1){120}} 
\multiput(10,120)(10,-10){10}{\circle*{2}} \put(10,120){\line(1,-1){90}} \put(10,120){\line(1,1){10}}
\multiput(10,100)(10,-10){9}{\circle*{2}} \put(10,100){\line(1,-1){80}} \put(10,100){\line(1,1){20}}
\multiput(10,80)(10,-10){8}{\circle*{2}} \put(10,80){\line(1,-1){70}} \put(10,80){\line(1,1){30}}
\multiput(10,60)(10,-10){6}{\circle*{2}} \put(10,60){\line(1,-1){50}} \put(10,60){\line(1,1){40}}
\multiput(10,40)(10,-10){3}{\circle*{2}} \put(10,40){\line(1,-1){20}} \put(10,40){\line(1,1){50}}
\multiput(10,20)(10,-10){2}{\circle*{2}} \put(10,20){\line(1,-1){10}} \put(10,20){\line(1,1){60}}

\put(20,10){\line(1,1){60}}
\put(50,20){\line(1,1){40}}
\put(60,10){\line(1,1){40}}
\put(80,10){\line(1,1){30}}

\put(5,120){\makebox(0,0)[r]{$[2]$}}
\put(5,100){\makebox(0,0)[r]{$[2,6]$}}
\put(5,80){\makebox(0,0)[r]{$[2,6,7]$}}
\put(5,60){\makebox(0,0)[r]{$[2,6,7,8]$}}
\put(5,40){\makebox(0,0)[r]{$[2,6,7,8,10]$}}
\put(5,20){\makebox(0,0)[r]{$[2,6,7,8,10,13]$}}

\put(140,5){\makebox(0,0)[t]{$\gamma_1$}}
\put(80,5){\makebox(0,0)[t]{$\gamma_2$}}
\put(60,5){\makebox(0,0)[t]{$\gamma_3$}}
\put(20,5){\makebox(0,0)[t]{$\gamma_4$}}
\end{picture}
\end{center}
\item
If $m=7$, $n=14$ and $\gamma=[2,6,7,8,10,13,14]$, then
$k=3$, $u(1)=1$, $u(2)=4$, $u(3)=5$ and
$\chi_0=\{1\}$,
$B_1=\{2\}$,
$\chi_1=\{3,4,5\}$,
$B_2=\{6,7,8\}$,
$\chi_2=\{9\}$,
$B_3=\{10\}$,
$\chi_3=\{11,12\}$
and
$B_4=\{13,14\}$.
$\gamma_1=[3,6,7,8,10,13,14]$, 
$\gamma_2=[2,6,7,9,10,13,14]$, 
$\gamma_3=[2,6,7,8,11,13,14]$ and 
$\gamma_4=[2,6,7,8,10,13]$.
Further,
$\xi(\gamma_1)=(5,0)$,
$\xi(\gamma_2)=(2,3)$,
$\xi(\gamma_3)=(1,4)$ and
$\xi(\gamma_4)=(-1,6)$.
Therefore, 
the Hasse diagram of $P\setminus\{\gamma\}$ is the following.
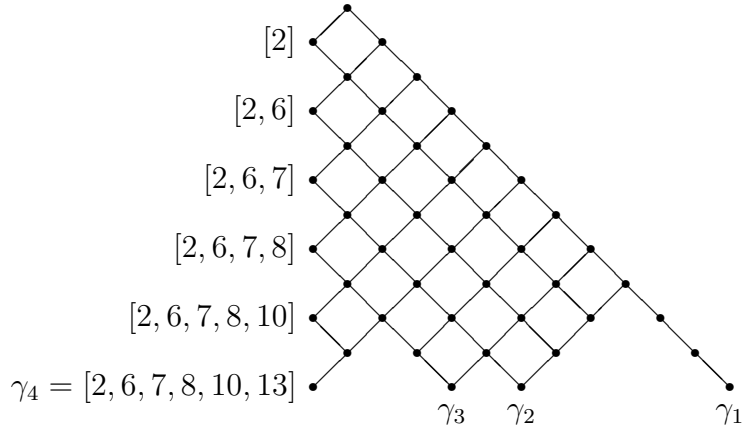
\begin{center}
\begin{picture}(170,130)(-30,0)

\multiput(20,120)(10,-10){12}{\circle*{2}} \put(20,120){\line(1,-1){110}} 
\multiput(10,110)(10,-10){9}{\circle*{2}} \put(10,110){\line(1,-1){80}} \put(10,110){\line(1,1){10}}
\multiput(10,90)(10,-10){8}{\circle*{2}} \put(10,90){\line(1,-1){70}} \put(10,90){\line(1,1){20}}
\multiput(10,70)(10,-10){7}{\circle*{2}} \put(10,70){\line(1,-1){60}} \put(10,70){\line(1,1){30}}
\multiput(10,50)(10,-10){5}{\circle*{2}} \put(10,50){\line(1,-1){40}} \put(10,50){\line(1,1){40}}
\multiput(10,30)(10,-10){2}{\circle*{2}} \put(10,30){\line(1,-1){10}} \put(10,30){\line(1,1){50}}
\multiput(10,10)(10,-10){1}{\circle*{2}} \put(10,10){\line(1,1){60}}

\put(40,20){\line(1,1){40}}
\put(50,10){\line(1,1){40}}
\put(70,10){\line(1,1){30}}

\put(5,110){\makebox(0,0)[r]{$[2]$}}
\put(5,90){\makebox(0,0)[r]{$[2,6]$}}
\put(5,70){\makebox(0,0)[r]{$[2,6,7]$}}
\put(5,50){\makebox(0,0)[r]{$[2,6,7,8]$}}
\put(5,30){\makebox(0,0)[r]{$[2,6,7,8,10]$}}
\put(5,10){\makebox(0,0)[r]{$\gamma_4=[2,6,7,8,10,13]$}}

\put(130,5){\makebox(0,0)[t]{$\gamma_1$}}
\put(70,5){\makebox(0,0)[t]{$\gamma_2$}}
\put(50,5){\makebox(0,0)[t]{$\gamma_3$}}

\end{picture}
\end{center}
\end{enumerate}
\end{example}

Now we state criteria of \gor\ property of the rings of 
absolute invariants.

First we consider the ring of absolute $\orth(m)$-invariants and 
give another proof of the result of Conca {\cite[Corollary 2.3]{con}}
(the case where $\gamma=[1,2,\ldots,m]$ is due to Goto \cite{got2}).
%
%
\begin{thm}
\label{thm:gor orth}
Suppose that $\gamma\neq[n-m+1,\ldots,n]$.
\begin{enumerate}
\item
If $b_m<n$, then
$K[Z_\gamma]^{\orth(m,-)}=K[Z_\gamma\transpose  Z_\gamma]$ is \gor\
if and only if 
$|B_i|=|\chi_{i-1}|$ for $i=2$, $3$, \ldots, $k$
and
$|\chi_k|\equiv1\pmod2$.
\item
If $b_m=n$, then
$K[Z_\gamma]^{\orth(m,-)}=K[Z_\gamma\transpose  Z_\gamma]$ is \gor\
if and only if 
$|B_i|=|\chi_{i-1}|$ for $i=2$, $3$, \ldots, $k$
and
$|B_{k+1}|=|\chi_k|-1$.
\end{enumerate}
\end{thm}
\begin{remark}\rm
We regard the null condition ``$|B_i|=|\chi_{i-1}|$ for $i=2$, \ldots, $k$''
valid if $k=1$.
\end{remark}
{\bf Proof of Theorem \ref{thm:gor orth}.}
Set $Q=\{[b_1,\ldots,b_i]\mid i=1$, \ldots, $m\}$.
Then $\ini K[Z_\gamma]^{\orth(m,-)}\simeq \DDDDD_K(\Gamma'(m\times n;\gamma), Q)$
by Theorem \ref{thm:ini is dhibi}.
Let $\tilde P$ be the poset defined as in the preceding paragraph of Theorem \ref{thm:gor cri}.

We first consider the case where $b_m<n$.
Then $\tilde P=P^+$ and $P^+$ is pure if and only if $|B_j|=|\chi_{j-1}|$
for $2\leq j\leq k$
by the argument before Example~\ref{ex:p setminus gamma}.
Further, since $\rank[\xi(y),\xi(y')]_{P'}\equiv 0\pmod 2$ for
any $y$, $y'\in Q^+\setminus\{\gamma\}$ with $y>y'$, 
\ref{item:rank even} of Theorem \ref{thm:gor cri} hold true if and only
if $\rank[\gamma,[b_1, \ldots, b_m]]\equiv\pmod 2$.
Since $\rank[\gamma,[b_1, \ldots, b_{m-1}]]=n-b_m+1=|\chi_k|+1$,
we see the result by Corollary \ref{cor:same hilb}
and Theorem \ref{thm:gor cri}.

Next we consider the case where $b_m=n$.
$\gamma\covered[b_1, \ldots, b_{m-1}]=\gamma_{k+1}$ in this case.
Thus \ref{item:rank even} of Theorem \ref{thm:gor cri} always holds true.
Further, $\tilde P$ is pure if and only if
$\height\xi(\gamma_1)=\cdots=\height\xi(\gamma_k)=\height\xi(\gamma_{k+1})+1$
in $P'$.
Therefore, 
by the argument before Example~\ref{ex:p setminus gamma},
we see the result by Corollary \ref{cor:same hilb}
and Theorem \ref{thm:gor cri}.
\qed

Note that if $\gamma=[n-m+1,\ldots, n]$, then we see that 
$K[Z_\gamma\transpose Z_\gamma]$ is isomorphic to a polynomial ring with
$\frac{1}{2}m(m+1)$ variables by
Lemma \ref{lem:lm of zgamma t} \ref{item:lm of zgamma t isom}.
In particular, $K[Z_\gamma\transpose Z_\gamma]$ is \gor.

Next we consider the ring of absolute $\sorth(m)$-invariants.
\begin{thm}\label{thm:gor sorth}
\begin{enumerate}
\item
If $b_m<n$, then
$K[Z_\gamma]^{\sorth(m,-)}=K[Z_\gamma\transpose  Z_\gamma,\Gamma(Z_\gamma)]$ is \gor\
if and only if 
$|B_i|=|\chi_{i-1}|$ for $i=2$, $3$, \ldots, $k$.
\item
If $b_m=n$ and $\gamma\neq[n-m+1,\ldots,n]$, then
$K[Z_\gamma]^{\sorth(m,-)}=K[Z_\gamma\transpose  Z_\gamma, \Gamma(Z_\gamma)]$ is \gor\
if and only if 
$|B_i|=|\chi_{i-1}|$ for $i=2$, $3$, \ldots, $k+1$.
\item
If $\gamma=[n-m+1,\ldots,n]$, then
$K[Z_\gamma]^{\sorth(m,-)}=K[Z_\gamma\transpose  Z_\gamma, \Gamma(Z_\gamma)]$ is \gor.
\end{enumerate}
\end{thm}
\begin{proof}
Set $Q=\{[b_1,\ldots, b_i]\mid 1\leq i\leq m-1\}$.
Then $K[Z_\gamma]^{\sorth(m,-)}\simeq \DDDDD_K(\Gamma'(m\times n;\gamma),Q)$
by Theorem \ref{thm:ini is dhibi}.
Let $\tilde P$ be the poset defined as in the preceding paragraph of 
Theorem \ref{thm:gor cri}.
Then $\tilde P=P^+$ and $\rank[y,y']_{\tilde P}\equiv 0\pmod 2$
for any $y$, $y'\in Q^+$ with $y<y'$,
i.e., \ref{item:rank even} of Theorem \ref{thm:gor cri} is always valid.
Therefore, 
if $\gamma\neq[n-m+1,\ldots, n]$, we see 
by Corollary \ref{cor:same hilb} and Theorem \ref{thm:gor cri},
that
$K[Z_\gamma]^{\sorth(m,-)}$ is \gor\ if and only if
$\height\xi(\gamma_1)=\cdots=\height\xi(\gamma_k)$ in $P'$
in the case where $b_m<n$
or 
$\height\xi(\gamma_1)=\cdots=\height\xi(\gamma_{k+1})$ in $P'$
in the case where $b_m=n$.
Further, we see that $K[Z_\gamma]^{\sorth(m,-)}$ is \gor\ if 
$\gamma=[n-m+1,\ldots,n]$.
Thus, the theorem follows
by the argument before Example~\ref{ex:p setminus gamma}.
\end{proof}
As a corollary, we see the following fact.

\begin{cor}
Suppose that $\gamma\neq[n-m+1,\ldots, n]$.
\begin{enumerate}
\item
If $b_m<n$ and $K[Z_\gamma]^{\orth(m,-)}$ is \gor,
then $K[Z_\gamma]^{\sorth(m,-)}$ is \gor.
\item
If $b_m=n$ and $K[Z_\gamma]^{\orth(m,-)}$ is \gor,
then $K[Z_\gamma]^{\sorth(m,-)}$ is not \gor.
\end{enumerate}
\end{cor}
%
%

%
%


\begin{thebibliography}{DEP2}
%
%
%
%
\bibitem[BC]{bc}
Bruns, W. and Conca, A.:
{\it
$F$-rationality of determinantal rings and their Rees rings.}
Michigan Math.\ J.\ {\bf 45} (1998) 291--299.
%
%
\bibitem[BH]{bh}
Bruns, W. and Herzog, J.:
``Cohen-Macaulay rings.''
Cambridge studies in advanced mathematics
{\bf 39} 
Cambridge University Press
(1993)
%
%
\bibitem[BV]{bv}
Bruns, W. and Vetter, U.:
``Determinantal Rings.''
Lecture Notes in Mathematics {\bf 1327} Springer (1988)
%
%
\bibitem[BRW]{brw}
Bruns, W. R\"omer, T. and Wiebe, A:
{\it
Initial algebras of determinantal rings, Cohen-Macaulay and Ulrich ideals.}
Michigan Math.\ J.  {\bf 53} (2005), 71--81.
%
%
\bibitem[Con]{con}
Conca, A:
{\it Divisor class group and the canonical class of determinantal
rings defined by ideals of minors of a symmetric matrix}
Arch.\ Math.\ {\bf 63} (1994), 216--224.
%
%
\bibitem[CHV]{chv}
Conca, A., Herzog, J.\ and Valla, G.:
{\it
Sagbi bases with applications to blow-up algebras.}
J.\ reine angew.\ Math.\ {\bf474} (1996), 113--138.
%
%
\bibitem[DEP1]{dep1}
DeConcini, C., Eisenbud, D. and Procesi, C.:
{\it Young Diagrams and Determinantal Varieties.}
Inv.\ Math.\ {\bf56} (1980), 129--165
%
%
\bibitem[DEP2]{dep2}
DeConcini, C., Eisenbud, D. and Procesi, C.:
``Hodge Algebras.''
Ast\'{e}risque
{\bf 91} (1982)
%
%
\bibitem[DP]{dp}
DeConcini, C. and Procesi, C.:
{\it A characteristic-free approach to invariant theory,}
Adv.\ Math.\ {\bf21} (1976), 330--354.
%
%
%
%
\bibitem[Got]{got2}
Goto, Shiro. "On the Gorensteinness of determinantal loci." Journal of Mathematics of Kyoto University 19.2 (1979): 371-374.
%
%
\bibitem[Hib]{hib}
Hibi, T.:
{\it Distributive lattices, affine smigroup rings and algebras 
with straightening laws.}
in ``Commutative Algebra and Combinatorics'' (M. Nagata and H. Matsumura, ed.),
Advanced Studies in Pure Math. {\bf11} North-Holland, Amsterdam (1987),
93--109.
%
%
\bibitem[Hoc]{hoc}
Hochster, M.:
{\it Rings of invariants of tori, Cohen-Macaulay rings generated
by monomials and polytopes.}
Ann. of Math. {\bf 96} (1972), 318-337
%
%
%
%
%
\bibitem[Miy]{miy1}
Miyazaki, M:
{\it 
A Note on the Invarints of the Group of Triangular Matrices
Acting on Generic Echelon Matrices,}
Bulletin of Kyoto University of Education Ser. B {\bf85} (1994), 57--62.
%
%
%
%
%
%
%
%
%
%
\bibitem[Ric]{ric}
Richman, David R.:{\it The fundamental theorems of vector invariants.}
Adv. Math. {\bf 73} (1989), 43--78.
%
%
\bibitem[Sta]{sta2}
Stanley, R. P.:
{\it Hilbert Functions of Graded Algebras.}
Adv. Math. {\bf 28} (1978), 57--83.
%
%
%
%
%
%
\end{thebibliography}
\end{document}